\documentclass[preprint,10pt]{elsarticle}

\usepackage{amsmath}
\usepackage{amssymb}
\usepackage{amsthm}

\usepackage[mathscr]{eucal}
\usepackage{graphicx}
\usepackage{float}
\restylefloat{figure}

\usepackage{fullpage}

\newtheorem{theorem}{Theorem}
\newtheorem{lemma}{Lemma}
\newtheorem{remark}{Remark}
\newtheorem{proposition}{Proposition}
\newtheorem{corollary}{Corollary}
\newtheorem{assumption}{Assumption}

\newtheorem{example}{Example}

\journal{Mathematics and Computers in Simulation}

\begin{document}

\begin{frontmatter}

\title{Noisy Prediction-Based Control Leading to Stability Switch}

\author[1]{E. Braverman$^3$}
\author[2]{ A. Rodkina}

\affiliation[1]{organization={University of Calgary}, city={Calgary}, postcode={T2N 1N4},  state={Alberta},  country={Canada}  }

\affiliation[2]{organization={The University of the West Indies, Mona Campus},  city={Kingston}, country={Jamaica}}




\begin{abstract}
Applying Prediction-Based Control (PBC) $x_{n+1}=(1-\alpha_n)f(x_n)+\alpha_n x_{n}$
with stochastically perturbed control coefficient $\alpha_n=\alpha+\ell \xi_{n+1}$, $n\in \mathbb N$, where $\xi$ are bounded identically distributed independent random variables, we globally stabilize the unique equilibrium $K$ of the equation
$
x_{n+1}=f(x_n)
$
in a certain domain.
In our results, the noisy control $\alpha+\ell \xi$ provides both local and global stability,
while the mean value  $\alpha$ of the control does not guarantee global stability, for example, the deterministic controlled system can have a stable two-cycle, and non-controlled map be chaotic. In the case of  unimodal $f$ with a negative Schwarzian derivative, we get sharp 
stability results generalizing Singer's famous statement `local stability implies global' to the case of the stochastic control.
New global stability results are also obtained in the deterministic settings for variable $\alpha_n$ and, generally, continuous but not differentiable at $K$ map $f$.
\end{abstract}


\end{frontmatter}


\section{Introduction}
\label{sec:Intr}

\footnotetext[3]{Corresponding author, e-mails maelena@ucalgary.ca; maelena@math.ucalgary.ca}

We consider the map
\begin{equation}
\label {eq:originalintr}
x_{n+1}=f(x_n), \quad x_0>0,
\end{equation}
where $f:[0, \infty)\to [0, \infty)$ is a continuous function  with one positive unstable equilibrium $K$,
$f(x)>x$ for $x\in (0, K)$,  $0<f(x)<x$ for $x\in (K, \infty)$.
Non-negativity of $x_n$ is assumed following a long tradition of population dynamics models, and the equilibrium $K$ of $f$ is unstable; moreover, $f$ can exhibit chaotic behaviour for maps such as Ricker, logistic, and others.

Various methods were applied to alleviate chaotic behaviour, 
some of them combined the current and the past values of $x_n$. 
In contrast to this approach, Prediction-Based Control (PBC) proposed by Ushio and Yamamoto \cite{uy99} computed the weighted average between the state variable $x$ and some iterate of the map $f^k(x)$
(a predicted, or a potential
 future value)
$
x_{n+1}=(1-\alpha)f^k(x_n)+\alpha x_{n},
$
which in the simplest case $k=1$ is 
\begin{equation}
\label{eq:PBCintrdet}
x_{n+1}=(1-\alpha)f(x_n)+\alpha x_{n}, \quad x_0>0, \quad \alpha \in [0,1).
\end{equation}
PBC was proved to be an efficient stabilization tool \cite{FL2010,polyak,sousa}. Moreover, some modification was considered recently in \cite{Chagas}.

While generally for parameter-based stabilization, increasing $\alpha$ does not lead to stability of the unique positive equilibrium point $K$, i.e. for a stabilizing $\beta_*$, instability can be observed for some $\alpha \in (\beta_*,1)$, there is a 
critical $\beta_*$, such that for any $\alpha \in (\beta_*,1)$, 
$K$ is a globally stable equilibrium of \eqref{eq:PBCintrdet}. 

Consider the case when $f$ is a three times differentiable unimodal function with two equilibrium points zero and $K$, a unique critical point $c\in (0,K)$ (maximum), $f'(0)>1$, $f''(x)<0$ for all $x \in (0,c)$ and a Schwarzian derivative
\begin{equation}
\label{eq:Schwarzian}
(Sf)(x) = \frac{f'''(x)}{f'(x)} - \frac{3}{2} \left( \frac{f''(x)}{f'(x)}  \right)^2,
\end{equation}
which is negative on $(0,\infty)$ excluding the unique critical point $c$. Under these conditions, local stability of $K$ implies global stability, and once $f'(K) < -1$, the sharp stabilizing constant $\displaystyle \alpha_0 := \frac{ -f'(K) - 1}{-f'(K)+1}$ determines the minimal value leading to stabilization \cite{FL2010}. The idea goes back to \cite{Singer}. As local stability guaranteed by the fact that the derivative of the right-hand side of \eqref{eq:PBCintrdet} at $K$ which is a weighted average of $f'(K)$ and one exceeds -1, is easily established, this leads to the fact that for any $\alpha \in (\alpha_0,1)$, $K$ is a globally stable equilibrium of \eqref{eq:PBCintrdet}, see \cite{FL2010}.
The main result of the paper is that this constant is no longer sharp in the stochastic case and can be improved when the control is perturbed by noise, which is rigorously justified for a symmetric continuous or discrete distribution. 

Our main goal is to  stabilize  globally the equilibrium $K$ applying PBC with stochastically perturbed variable control coefficient $\alpha_n=\alpha+\ell \xi_{n+1}$, $n\in \mathbb N$,  where $\xi$ are bounded identically distributed independent random variables,
\begin{equation}
\label {eq:PBCintr}
x_{n+1}=(1-\alpha- \ell \xi_{n+1})f(x_n)+(\alpha+\ell\xi_{n+1})x_{n}, \quad x_0>0.
\end{equation}
Once $f$ satisfies some smoothness criterion being at least one-sided Lipschitz continuous at $K$, such control always exists even with $\ell=0$. However, our purpose is to find the smallest  possible value of the parameter $\alpha$ and the range for the noise level $\ell$ which provide global stability of $K$ with probability one (almost surely).
This is a two-parameter $(\alpha,\ell)$-problem, also dependent on the choice of the distribution for $\xi$.
A bound is constructed for a couple $(\alpha,\ell)$ guaranteeing convergence of a solution to 
the unique positive equilibrium. For a map with negative Schwarzian derivative~\eqref{eq:Schwarzian}, we get a sharp stabilization criterion which is, unlike most stochastic results, global. It also clearly illustrates that the range of $\alpha$ includes smaller control values than in the deterministic case \cite{FL2010,Singer}.

Let us notice that local stabilization of an unstable equilibrium is possible with noise only.
However, without any other control, for a chaotic map the attracting neighbourhood of $K$ can be very small, less than $10^{-9}$~\cite{BR2019}.
Most relevant results on local stability of the stochastic difference equation compared to the present paper can be found in \cite{BR2019,Medv}, where unbounded noise $\xi$ and continuously differentiable $f$ were considered (see also \cite{BR0, BR1, BKR2020}). This approach is due to Khasminskii \cite{Kh} and  H. Kesten \cite{K};
it was used  later in many publications.   Our main theorems are based on the Kolmogorov's Law of Large Numbers which is used to compute the values of $\ell$ and $\alpha$ ensuring local stability with any given probability, 
and on a corollary of the Borell-Cantelli Lemma, guaranteeing that 
a solution eventually enters however small, but prescribed in advance, 
neighbourhood of $K$, leading to global stability with the probability one.

Stabilizing effect of noise attracted a lot of attention recently due to its significant role in sustaining 
healthy neuronal activities and avoid sustainable oscillations \cite{Rich}. In contrast to \cite{Rich}, we consider bounded, not Gaussian noise, which is assumed to be more realistic in biological and health-related systems.

The main result of the paper proves global asymptotic stability of the solution to \eqref{eq:PBCintr} under the assumption that there is local stochastic asymptotic stability of $K$ with probability not less than $1-\gamma$ in the interval  $(K-\delta, K+\delta)$ for the initial value, where $\delta=\delta(\gamma)$ depends on $\gamma\in (0, 1)$ chosen arbitrarily. The control parameters $\alpha$ and $\ell$ have also to satisfy  $\alpha+\ell>\beta_*$,  $\alpha<\beta_*$, where $\beta_*$ is  a control level ensuring  global asymptotic stability for deterministic equation \eqref{eq:PBCintrdet}.
The control $\alpha+\ell \xi$ provides both local and global stability, 
while the mean value  $\alpha$  does not necessarily guarantee global (and maybe even local) stability of $K$.

Continuity of $f$ and the fact that $(f(x)-K)(x-K)<0$ on $(0,K) \cup (K,\infty)$ allow us to deal with a smaller interval $[f^2_m, f_m]$ around $K$, where  $f_m$ is a maximum of $f$ on $[0,K]$, and $f^2_m$ is a minimum of $f$ on $[K, f_m]$. Then,
the solution to \eqref{eq:PBCintr}  reaches $[f^2_m, f_m]$ in a finite number $S_0$ of steps and remains there. Keeping a control level greater than  $\beta_*$ for some prescribed finite number of steps allows a solution $x$ to get from $[f^2_m, f_m]$ into the  initial interval $(K-\delta, K+\delta)$, from where $x$ converges to $K$.  Applying Borel-Cantelli  lemma (see Lemma~\ref{lem:topor}), we conclude that there exists a random moment $\mathcal N$ 
such that the required control intensity occurs for a certain prescribed number of steps in a row, which implies global stability. 

Generally speaking, the function $f$ is not assumed to be continuously differentiable at $K$, which allows to distinguish between the left-side Lipschitz-type constant $L^-$ and the right-side $L^+$ at $K$, and find a control parameter for local, as well as global stability, using $L^-$ and $ L^+$, instead of their maximum.  The simplest result about local stability in the deterministic setting is obtained when 
$\displaystyle \alpha \in \left(\underline \alpha_0, 1\right)$, $\underline \alpha_0:=\frac{L^+L^--1}{(L^++1)(L^-+1)}$, which, to the best of our knowledge,  is a new result. For $L^-=L^+=L$, a well-known condition of $\alpha> (L-1)/(L+1)$ follows, see e.g. \cite[Remark 1]{BL2012}.

When local stability is due to the stochastic control, the proof applies the Kolmogorov's Law of Large Numbers (see Lemma~\ref{lem:Kolm}), once
\begin{equation}
\label{cond:locintr}
\mathbb E \ln \left([L^--(\alpha+\ell \xi))(L^-+1)] [L^+-(\alpha+\ell \xi))(L^++1)]  \right)<0
\end{equation}
is satisfied. For some distributions of $\xi$, we show that when $\ell>0$, condition \eqref{cond:locintr} gives smaller value of $\alpha$ than $\underline \alpha_0$. Notice that if  $f$ is a unimodal function satisfying the conditions of \cite{FL2010} elaborated above, any parameter $\alpha> \alpha_0:=\frac{-f'(K)-1}{-f'(K)+1}$ provides local, as well as global stability  of deterministic equation  \eqref{eq:PBCintrdet}. We calculate the values of $\ell$ for certain types of noise $\xi$ such that the global stability for \eqref{eq:PBCintr} holds with some 
$\alpha < \alpha_0$.

Even though there is always $\beta_*$ which guarantees  global stability for \eqref{eq:PBCintrdet}, we are interested in the smallest possible one. In many population biology models, like chaotic Ricker and logistic, using just the left-hand-side  $L^-$ and the right-hand-side $L^+$ global Lipschitz-type constants  gives much bigger value than  
$\alpha_0:=\frac{-f'(K)-1}{-f'(K)+1}$, as $\frac{L^+L^--1}{(L^++1)(L^-+1)}> \alpha_0$.  It appeared advantageous  to split $[f^2_m, K]$ into several subintervals  and pair each one with  its image by the map 
$(1- \alpha_0)f(x)- \alpha_0$, calculating $L^-_i$ and $L^+_i$, Lipschitz constants for the left and the right intervals.
The point is that, while the left constants are very large, low right constants alleviate for them.
For Ricker and bobwhite quail \cite{mb} models, $L^-_i$ and $L^+_i$ compensate each others: when the left one is getting bigger, the right one is getting smaller, so the expression $\frac{L^-_iL^+_i-1}{(L^+_i+1)(L^-_i+1)}$ which can be used  for calculating $\beta_*$, remains less than a control value $\alpha_0$ computed without this splitting. 

Our approach to the proof of global stability for the deterministic equation is based on 
the results of \cite{Gull,Cull2,Shark}, see also \cite{ElaydiSacker,Coppel}.
As discovered in  \cite{Coppel}, a unique equilibrium of $f$ is globally stable if and only if $f^2$ has no two-cycles. We cite this statement as in \cite[Corollary C.4]{Elaydi}.

\begin{lemma}
\label{lem:Gull}
Let $g:[a,b]\to [a,b]$  be continuous, then its fixed point $x^*$
is globally asymptotically stable  relative to $[a,b]$ if and only if $g^2( x)>x$, $x<x^*$ and  $g^2(x)<x$, $x>x^*$
, for all $x\in (a, b)\setminus \{x^*\}$, and either $g(a)<b$ or $g(b)>a$.
\end{lemma}

When the conditions of \cite{FL2010} are satisfied, local stability implies global stability
in the deterministic case. We extend this sharp result to the stochastic case: some version of local condition \eqref{cond:locintr} implies global stability with probability 1. While the main results of the paper refer to global stability of stochastically perturbed maps, there are new findings for deterministic equations with variable control, or for continuous non-smooth$f$. 

In  \cite{BR2022},  PBC was used to stabilize simultaneously multiple   equilibrium points of \eqref{eq:originalintr}. It was supposed that $f(x)-x$ changes its sign at each $K_j$, $j=0, \dots j_0$, and at each $K_{2i+1}$,  $f$ satisfies a one-side Lipschitz condition. The control was defined, based on the minimum $L$ of  the  left and right-side Lipschitz constants, whenever available.
There was a total of $j_0\ge 4$ equilibrium points $K_j$,   and every second one, $K_{2i+1}$, was stabilized.  The minimum number of stabilized equilibriums was 2, when $j_0=4$, so the case of the unique  positive equilibrium was not considered in  \cite{BR2022}.  
Compared to \cite{BR2022},  the present paper has the following common features: sharp results are achieved by careful computation of the minimal control constant allowing to avoid a two-cycle; also, common tools are used, in particular, the proofs in the case when the control is stochastically perturbed are based on the Borel-Cantelli Lemma.
However, the models are different: roughly speaking, \cite{BR2022} is focused on pulse stabilization with PBC 
at each $2^k$th step, while we focus on the original map and classical PBC.  While conditions in \cite{BR2022} are not easy to verify, we obtain sharp results for smooth unimodal functions with a negative Schwarzian derivative when global stability can be established by checking easily verifiable local conditions, and our method of finding the best control in deterministic settings is different in the present paper from \cite{BR2022}. 

The rest of the paper is structured as follows. In Section~\ref{sec:defas}, we  formulate properties of the noises $\xi_n$ and state two important results applied in the proof of  the main theorems: the Kolmogorov's Law of Large Numbers and a corollary of the Borel-Cantelli Lemma. 
We also define properties of auxiliary functions that are later used in the paper.  
Section~\ref{sec:globdet}  discusses  global stability for deterministic equation \eqref{eq:PBCintrdet}  when control parameter is  either constant or  variable. 
In Section~\ref{sec:GLstoch} we establish conditions when local stability implies  the global stability when 
the control is stochastically perturbed as in \eqref{eq:PBCintr}: in Section~\ref{subsec:detloc} local stability holds for the deterministic value of  control  while in Section \ref{sec:genglobloc} stochastic local stability was obtained by application of the   Kolmogorov's Law of Large Numbers, which is the main result of the paper. 
In the case when $f$ is unimodal with a negative Schwarzian derivative \cite{FL2010}, in Section~\ref{subsec:Zinger} we show that stochastic perturbation of the control can improve the sharp deterministic constant for the average control.
In Section~\ref{subsec:parglobstab}, we calculate the  control parameter which provides global stability based on the left and the right global Lipschtiz constants and generalize this to several intervals with different Lipschitz constants.
Section~\ref{sec:exsim} contains examples and simulations which illustrate our results.
All proofs are deferred to the Appendix. 

\section{Assumptions and Auxiliary Statements}
\label{sec:defas}

Denote by
$\left[x\right] $ the largest integer not exceeding $x$, 
${\mathbb N}_0:=\mathbb N\cup \{0\}$, ``s.t'' stands for ``such that''.

\subsection{Assumptions on the noise}

Introduce a complete filtered probability space $(\Omega, {\mathcal{F}}$, $\{{\mathcal{F}}_n\}_{n \in 
\mathbb N}, {\mathbb P})$, where the filtration $(\mathcal{F}_n)_{n \in \mathbb{N}}$ is naturally generated by 
the sequence of independent identically distributed random variables $(\xi_n)_{n\in\mathbb{N}}$, i.e. 
$\mathcal{F}_{n} = \sigma \left\{\xi_{1},  \dots, \xi_{n}\right\}$.
The standard abbreviation ``a.s.'' is used for either ``almost sure" or ``almost surely" 
with respect to the probability measure $\mathbb P$, and 
``i.i.d.'' for  ``independent identically distributed'', to describe random variables. 
For a detailed introduction to stochastic concepts and notations, we refer the reader to \cite{Shiryaev96}. 

In this paper we consider bounded noises and control perturbations, which is a natural assumption in  population dynamics. 
\begin{assumption}
\label{as:noise}
$(\xi_n)_{n\in \mathbb N}$ is a sequence of i.i.d. random variables such that $|\xi_n|\le 1$,  $\forall n\in \mathbb N$, and, for each $\varepsilon>0$, $\mathbb P \{\xi\in (1-\varepsilon, \, 1] \}>0$.
\end{assumption}

The following lemma was proved in \cite{BRAllee} and is a corollary of the Borel-Cantelli Lemma.
\begin{lemma}
\label {lem:topor} 
Let sequence $(\xi_n)_{n\in \mathbb N}$ satisfy Assumption \ref{as:noise}. Then, for each nonrandom $S\in \mathbb N$, $\varepsilon\in (0, 1)$ and a random moment $\mathcal M$ we have 
$$
\mathbb P\{ \mbox{there exists a random moment}  \quad  \mathcal N>\mathcal M:\xi_{\mathcal N+i} \in (1-\varepsilon, 1), \,\,  i=0, 1, \dots, S \}=1.
$$
\end{lemma}

For simulations in the present paper, we consider discrete, as well as continuous random variables $\xi_n$ with a symmetric distribution.
 As an example of discrete distribution, we use Bernoulli random variable $\xi$, which takes the values of 1 and -1 with probability 1/2 each, has the zero mean  and the second moment  $\mu_2=1$. As an example of continuous random variable we use  continuous uniformly distributed on $[-1, 1]$ random variable $\xi$, which has the mean zero and the second moment  $\mu_2=1/3$.

The Kolmogorov Law of Large Numbers is cited below, see Shiryaev~\cite[P. 391]{Shiryaev96}. 
\begin{lemma}
\label{lem:Kolm}
  Let $(v_{n})_{n\in\ \mathbb N}$ be a sequence of i.i.d. random
  variables with $\mu:=\mathbb E v_n $, $\mathbb  E|v_n|<\infty$, $n\in \mathbb N$. Then, a.s.,~~ $\displaystyle \frac{1}{n}\sum_{k=1}^n v_k \rightarrow \mu$ as $n \to \infty$.
\end{lemma}

\begin{corollary}
\label{cor:Kolm}
Under the assumptions of Lemma~\ref{lem:Kolm}, for each $\varepsilon>0$ there exists a random $\mathcal N(\varepsilon)$ such that\
\begin{equation}
\label{est:Kolmv}
(\mu-\varepsilon)n\le \sum_{k=1}^n v_k\le  (\mu+\varepsilon)n, \, \,\,\mbox{for}\, \, \,n\ge \mathcal N(\varepsilon), \,\,\, \mbox{a.s.}
\end{equation}
Also, for each $\gamma\in (0, 1)$, there exist a nonrandom $\rm N=\rm N(\gamma, \varepsilon)$ and  $\Omega_\gamma \subset \Omega$ with $ \mathbb P(\Omega_\gamma)>1-\gamma$, such that \eqref{est:Kolmv}  holds on $\Omega_\gamma$ for $n\ge N$.
\end{corollary}

\subsection{Assumptions on $f$ and properties of  auxiliary functions}
\begin{assumption}
\label{as:Lglob1}
Let $f: [0, \infty)\to [0, \infty)$ be a continuous function with one positive locally unstable equilibrium $K>0$, $f(x)>x$ for $x\in (0, K)$,  $0<f(x)<x$ for $x\in (K, \infty)$, $f(0)=0$.
\end{assumption}
For $f$ satisfying Assumption~\ref{as:Lglob1}, we introduce the values
\begin{equation*}
\begin{split}
& x_{\max} \,\, \mbox{is the largest point of maximum of $f$ on $[0, K]$}, \,\, f_m:=\max\{f(x), \, x\in [0, K] \}=f(x_{\max}), \\
&
f^2_m:=\inf\{f(x), x\in (K, f_m) \}.
\end{split}
\end{equation*}

The constant $f_m$ is  well defined for any continuous $f$ and, since our purpose is to stabilize  the unstable equilibrium $K$, we have $0<f^2_m<K<f_m$. 

\begin{assumption}
\label{as:Lglob2}
The function $f$ satisfies Assumption \ref{as:Lglob1}, and for some $L_-\ge L_+>1$,
\begin{equation}
\label {cond:Lglob2}
\begin{split}
& f(x)-K\le L^- (K-x), \,\,  \mbox{\rm if} \,\,  x\in (x_m, K), \quad K-f(x)\le L^+ (x-K), \,\,  \mbox{\rm if} \,\, x\in (K, f_m).
\end{split}
\end{equation} 
\end{assumption}
Everywhere in the paper we assume $L_-\ge L_+>1$. 
Assumptions~\ref{as:Lglob1}-\ref{as:Lglob2} are  satisfied for many common population dynamics maps, such as unstable Ricker, logistic and Beverton-Holt models, where  $L_-\ge L_+>1$. Because of this, we concentrate on functions $f$ which are steeper to the left than to the right of $K$. The case $L^+ \ge L^->1$ is similar, so we omit discussing it. 

Define 
\begin{equation}
\label{def:Gbeta}
G(\beta, x):=(1-\beta)f(x)+\beta x, ~~ \beta\in [0, 1), ~~ x\geq 0.
\end{equation}
The properties of $G$, which are widely used in this paper, are stated in the next lemma, partially they were
justified in \cite[Lemma 2.2]{BR2022}. 
 \begin{lemma}
\label{lem:Gv}
Let  $f$ satisfy Assumption~\ref{as:Lglob1} and $G$ be defined as in \eqref {def:Gbeta}. Then
\begin{enumerate}
\item [(i)] $G(1,x)=x$, $G(0,x)=f(x)$, $x\in [0, \infty)$, $G:[0,1]\times \mathbb R\to  \mathbb R$ is a continuous function.
\item [(ii)] $f(x)>G(b, x)>G(a, x)>x$, if  $1>a>b>0$ and $x<K$, while 
$f(x)<G(b, x)<G(a, x)<x$, if  $1>a>b>0$ and $K<x$.
\item [(iii)] For $\beta\in(\beta_0, 1)\subset (0, 1)$, and  $\hat \beta:=\frac{\beta-\beta_0}{1-\beta_0}$, we have $
G(\beta, x)=(1-\hat \beta)G(\beta_0, x)+\hat \beta x.$
\item [(iv)]  $G(\beta,  \cdot): \bigl[f^2_m, \, f_m\bigr] \to \bigl[f^2_m , \, f_m\bigr]$ for all $\beta\in [0, 1)$.
\item [(v)] If $x_{G\max}$  is the largest point of maximum of $G(\beta, x)$ on $[0, K]$, $\beta\in [0, 1)$ then $x_{\max}\le x_{G\max}$.
\end{enumerate}
\end{lemma}
Set, for $\beta\in (0, 1)$, $L^{\pm}$  from  \eqref{cond:Lglob2} and  $L:=\max\{L^+, L^-\}$,
 \begin{equation}
\label{def:mathcalL_+12}
\mathcal L^{\pm}(\beta):=(1-\beta) L^{\pm}-\beta, \quad \mathcal L(\beta):=(1-\beta) L-\beta.
\end{equation} 
Define also
\begin{equation}
\label {def:Psi}
\Psi(u,v):=\frac{uv-1}{(u+1)(v+1)}=1-\frac 1{u+1}-\frac 1{v+1},\quad (u,v)\in (-1, \infty)\times (-1, \infty).
\end{equation} 
The next lemma states some useful properties of functions $\mathcal L^{\pm}(\cdot)$ and $\Psi(\cdot, \cdot)$. 

\begin{lemma}
\label{lem:mathcalL}
Let Assumptions \ref{as:Lglob1}-\ref{as:Lglob2} hold, and $L^->L^+$. Using  notations from \eqref{def:mathcalL_+12} and \eqref{def:Psi}, we have
\begin{enumerate}
\item [(i)] The functions $ \mathcal L^+(\beta)$ and  $\mathcal L^-(\beta)$ are monotone decreasing  for $\beta\in [0, 1]$.
\item [(ii)] $ \mathcal L^+(\beta)< \mathcal L^-(\beta)$  for $\beta\in [0, 1)$, \, $ \mathcal L^{\pm}(1)=-1$, $ \mathcal L^{\pm}\left( \frac{L^{\pm}}{L^{\pm}+1} \right)=0$, \,$ \mathcal L^{\pm}\left( \frac{L^{\pm}-1}{L^{\pm}+1} \right)=1$. 
\item [(iii)] $\mathcal L^+(\beta)\in [0,1], \quad \mathcal L^-(\beta)\in [1, \infty)$ when $\beta\in \left[\frac{L^+-1}{L^++1}, \min\left\{ \frac{L^--1}{L^-+1}, \frac{L^+}{L^++1} \right\}\right]$. 
\item [(iv)] $ \mathcal L^-(\beta)\le 1$ \, if \, $\beta\ge \frac{L^--1}{L^-+1}$, \, $ \mathcal L^{\pm}(\beta)\le 0$ \,  if\,  $\beta\ge \frac{L^{\pm}}{L^{\pm}+1}$.
\item [(v)] The function $\Psi(\cdot, \cdot)$ increases in each argument, does not exceed one  and  is positive for $uv>1$.
\item [(vi)]  The quadratic equation $\mathcal L^+(\beta)\mathcal L^-(\beta)-1=0 $  in $\beta$ has a positive discriminant, its smallest solution
\begin{equation}
\label {def:beta0}
\underline \beta_0=\Psi(L^- , L^+)=
\frac{L^- L^+-1}{(L^-+1)(L^++1)},
\end{equation}
while the largest solution is equal to 1. Moreover, for $L^- L^+>1$, we have $\underline \beta_0\in (0, 1)$,
$\underline \beta_0\in \left[\frac{L^+-1}{L^++1}, \min\left\{ \frac{L^--1}{L^-+1}, \frac{L^+}{L^++1} \right\}\right],$ and $\underline \beta_0=\frac{L^+-1}{L^++1}$ if $L^+=L^-$.
\item [(vii)]
The inequality  $\mathcal L^+(\beta)\mathcal L^-(\beta)<1 $ holds when $\beta\in (\underline \beta_0, 1)$. Moreover,
 $\mathcal L^+(\beta)\mathcal L^-(\beta)<1 $ if and only if $\Psi(L^- , L^+)<\beta.$
  \end{enumerate}
\end{lemma}
In terms of $G$, defined as in \eqref{def:Gbeta}, equations \eqref{eq:PBCintrdet} and \eqref{eq:PBCintr} can be written as 
\begin{equation}
\label {eq:PBCGstoch}
x_{n+1}=G(\alpha_n, x_n)=(1-\alpha_n)f(x_n)+\alpha_n x,  \quad x_0>0,\quad n\in \mathbb N,
\end{equation}
where $\alpha_n\in [0,1]$ is a variable parameter, which can be deterministic or stochastic in the form $\alpha_n=\alpha+\ell \xi_{n+1}$. 

Next lemma shows that $\bigl[f^2_m, \, f_m\bigr] $ is a trap which can be reached after a finite explicitly computed number of steps.

\begin{lemma}
\label{lem:aux21}
Let Assumption \ref{as:Lglob1} hold, $\beta_*\in (0, 1)$,  $\beta^*\in (\beta_*, 1)$ and $\beta_*\le \alpha_n\le \beta^*$. Let $x$ be a solution to \eqref{eq:PBCGstoch}  with $x_0>0$.  Then there exists a finite number $S_0=S_0(\beta_*, \beta^*, x_0)\in \mathbb N$  such that $x_n\in \bigl[f^2_m, \, f_m\bigr] $, for $n\ge S_0$.
\end{lemma}
 
\section{Global stability of the deterministic equation}
\label{sec:globdet}

In this section we discuss global stability of deterministic  equations \eqref{eq:PBCintrdet} and \eqref{eq:PBCGstoch}. Note that  global stability of $K$ for $\alpha \in \left( \frac{L^-}{L^-+1} , 1\right)$ was proved in 
\cite{BL2012}, and generalized to  variable $\alpha$ in \cite{BKR2016}. 
The results of this section are based on Lemma~\ref{lem:Gull} and \cite[Theorem 3]{Cull2}.  

\begin{lemma}
\label{lem:Gullbeta}
Let Assumption \ref{as:Lglob1} hold, $L^->L^+$ and $G(\beta, x)$  be defined as in \eqref{def:Gbeta} and
\begin{equation}
\label {def:beta*Gull}
\underline \beta_*:= \inf \mathcal S, \mbox{~~where~~} \mathcal S:=\{\beta\in (0, 1):   G^2(\beta, x)<x, ~x\in (f^2_m, K), ~~  G^2(\beta, x)>x,~x\in (K, f_m)\}.
\end{equation}
Then $\mathcal S$ is non-empty, and for any $\alpha \in (\underline \beta_*,1)$, 
any solution $x$ to \eqref{eq:PBCintrdet} with $x_0>0$ satisfies $\lim\limits_{n\to \infty}x_n=K$.
\end{lemma}

\begin{remark}
\label{rem:IFFGull}
Lemma~\ref{lem:Gull} and the definition of $\mathcal S$ in \eqref{def:beta*Gull}  yield  that the  parameter $\underline \beta_*$   is sharp  for deterministic equation \eqref{eq:PBCintrdet} to guarantee global asymptotic stability.  In general,  finding the best control $\underline \beta_*$ is not  trivial. 
However, stochastic perturbations still can decrease the mean value $ \underline \beta_*$  of the stabilizing control, see 
Sections~\ref{subsec:detloc}, \ref{subsec:Zinger} and relevant examples.

Note that definition \eqref{def:beta*Gull} does not exclude existence of the point $x_*\in [f^2_m, f_m]$ such that $G^2( \underline \beta_*, x_*)=x_*$, so when $\alpha=\underline \beta_*$  the solution to \eqref{eq:PBCintrdet} is not  globally asymptotically stable and there is a two-cycle $(x_*, G( \underline \beta_*, x_*))$. For each $x\in [f^2_m, K]$, we have $G^2( \underline \beta_*, x)\ge x$, while  $x\in [K, f_m]$ implies $G^2( \underline \beta_*, x)\le x$.
 \end{remark}
The following theorem is the main result of this section.  \begin{theorem}
   \label{thm:globalvarcontr} 
   Let Assumption \ref{as:Lglob2}  hold,  $\underline \beta_*$ be defined as in \eqref{def:beta*Gull}, $\beta_*\in (\underline \beta_*, 1)$,
    $\beta^*\in (\beta_*, 1)$ and   $\alpha_n\in [\beta_*, \beta^*]$.  
    Then
   \begin{enumerate}
    \item [(i)] The solution to \eqref{eq:PBCGstoch}  with  any $x_0>0$ converges to $K$.
    \item [(ii)] For any  $x_0>0$ and $\delta>0$ there is a finite number of steps $S_1=S_1(x_0, \beta_*, \delta)$ s.t. $x_n\in (K-\delta, K+\delta)$ for $n\ge S_1$.
 \end{enumerate}	    
   \end{theorem} 
   The proof of Theorem \ref{thm:globalvarcontr}  is based on Lemma~\ref{lem:Gull}  and  Lemma~\ref{lem:auxback} below.
   
   \begin{lemma}
   \label{lem:auxback}
    Let  Assumption  \ref{as:Lglob2}  hold,  $\underline \beta_*$ be defined as in \eqref{def:beta*Gull}, $\beta_*\in (\underline \beta_*, 1)$,  and  $\beta^*\in (\beta_*, 1)$.  Let  $x$ be a solution to \eqref{eq:PBCGstoch}   with   $\alpha_n\in [\beta_*, \, \beta^*]$. 
 If  $x_n<K$  for some $n\in \mathbb N$ and $s>n$ then $x_{s}>x_n$,  while if $x_n>K$ for some $n \in \mathbb N$  and $s>n$ then $x_{s}<x_n$. 
   \end{lemma}

\section{Global stability induced by noise}
\label{sec:GLstoch}

Now we proceed to the case when the control is stochastically perturbed as in \eqref{eq:PBCintr} and  establish global stability conditions. Let $\underline \beta_*$ be  defined as in \eqref{def:beta*Gull}, so that  deterministic equation  \eqref{eq:PBCintrdet} is  globally stable for any $\alpha>\underline \beta_*$.
We concentrate on the case when deterministic equation  \eqref{eq:PBCintrdet} is locally stable  for  $ \alpha>\alpha_0$, where $\alpha_0<\underline\beta_*$.
We show, under some restrictions, that introduction of noise into the control allows to get the stability result for  \eqref{eq:PBCintr}  for an intermediate value of the parameter $\alpha_0<\alpha <\underline\beta_*$ and a corresponding noise level $\ell>0$. In other words, stochastically perturbed control $\alpha+\ell \xi$ provides both local and global stability, while control with the mean value  $\alpha$  does not lead to global stability of $K$. 
 
 The shape of function $f$ guarantees that any solution gets into the interval $[f^2_m, f_m]$ after a finite number of steps $S_1$, which depends on the initial value $x_0$.  This interval is the first trap: a solution necessarily gets into this interval and stays there forever.  If the deterministic control satisfies $\beta>\underline\beta_*$, there is  a  finite number of steps $S_2$ after which a solution of  \eqref{eq:PBCintrdet}  gets into the second trap $(K-\delta_0, K+\delta_0)$, from where it converges to $K$. Thus, if  $\alpha+\ell>\underline\beta_*$,  by applying the Borel-Cantelli Lemma~\ref{lem:topor}, we conclude that starting from some random moment $\mathcal N$, control $\alpha+\ell\xi_n$ remains  greater than $\underline\beta_*$ for at least $S_2$ number of steps in a row. This allows a solution  of \eqref{eq:PBCintr} to get into the second trap. 
 Therefore $x_n\in (K-\delta_0, K+\delta_0)$ for $n>S_1+\mathcal N+S_2$, and $\displaystyle \lim_{n\to \infty}x_n=K$.
 We have to stress that, in stochastic settings, even though local stability is established with any given probability $1-\gamma$, $\gamma\in (0, 1)$,  on respected initial values interval,  global attractivity holds a.s. 

In Section~\ref{subsec:detloc},  we state global stability when local stability is provided by the deterministic control  $\alpha_0$. In this case, the coefficients of the stochastically perturbed control $\alpha+\ell\xi$ are easily calculated, independently of the distribution of the noise $\xi$, see \eqref{cond:ael} below: the average part $\alpha$ which provides  the global control satisfies $\displaystyle \alpha > \frac{\alpha_0+\underline\beta_*}2$.
By application of Lemma~\ref{lem:Kolm}, in Section~\ref{sec:genglobloc} we prove a theorem on global stability when local stability is ensured by a stochastically perturbed control $\bar \alpha_0$.  Note that $\bar \alpha_0\le \alpha_0$, where $ \alpha_0$ provides deterministic local stability.
In this case it is not so easy to estimate the average part $\alpha$ of  the global control, it just should satisfy local condition \eqref{cond:lambda0} and some more restrictions.
In Section \ref{subsec:Zinger}, we consider a unimodal function $f$ with a negative Schwarzian derivative.
For such $f$ in deterministic setting, local stability with a parameter $\alpha_0$ implies the global one \cite{FL2010,Singer}. However, with noise, the minimum deterministic constant $\alpha_0$ providing stability is no longer sharp in the sense that the control $\alpha+\ell \xi$ stabilizes for some $\alpha<\alpha_0$ and $\ell>0$, which is illustrated in the cases of Bernoulli (taking the values of $\pm 1$ with the probability of 0.5 each) and continuous uniformly distributed on $[-1,1]$ types of noise.

\subsection{Deterministic local stability implies global stability}
\label{subsec:detloc}

\begin{theorem}
 \label{thm:locgldet}
 Let  Assumptions  \ref{as:noise},  \ref{as:Lglob1} hold, the value of  $\underline \beta_*$ be  defined as in \eqref{def:beta*Gull}, $\delta_0\in (0, K)$, and $\alpha_0\in (0, 1)$. Let a solution $x$ to \eqref{eq:PBCintrdet} with any $\alpha>\alpha_0$ and $x_0\in (K-\delta_0, K+\delta_0)$  satisfy $\displaystyle \lim_{n\to \infty} x_n=K$, and
 \begin{equation}
 \label{cond:ael}
 \alpha\in \left(\frac{\alpha_0+\underline\beta_*}2, \underline\beta_* \right), \quad \ell \in (\underline\beta_*-\alpha, \, \min\{\alpha-\alpha_0, 1-\alpha\}).
 \end{equation}
 Then a solution $x$ to \eqref{eq:PBCGstoch} with $(\alpha, \ell)$ from \eqref{cond:ael} and any $x_0>0$  satisfies $\displaystyle \lim_{n\to \infty} x_n=K$ a.s.
  \end{theorem}
Note that $\alpha_0$ can be found as $\alpha_0=\frac{-f'(K)-1}{-f'(K)+1}$ if $f$ is differentiable at $K$, 
and as $\alpha_0=\Psi(\tilde L^-, \tilde L^+)$ if it is not. 
The latter statement is confirmed in the next Lemma~\ref{lem:locdet}.

Let, instead of Assumption~\ref{as:Lglob2},  local Lipschitz-type conditions hold:
\begin{equation}
\label {cond:loc2}
f(x)-K\le \tilde L^-[K-x] \quad \mbox{if} \,\,   x\in [K-\theta, K], \quad K-f(x)\le \tilde L^+[x-K] \quad\mbox{if} \,\,  x\in [K, K+\theta].
\end{equation}
\begin{lemma}
\label{lem:locdet}
Let conditions  \eqref{cond:loc2} hold,  $\Psi$ be defined in \eqref{def:Psi}, $\alpha>\Psi(\tilde L^+, \tilde L^-)$ and $\tilde{\mathcal L}^-(\alpha)$  be as in \eqref{def:mathcalL_+12}. Then $\lim\limits_{n\to \infty}x_n=K$ for each solution $x$ to \eqref{eq:PBCintrdet} for any $x_0\in (K-\theta/ \tilde{\mathcal L}^-(\alpha), K+\theta)$.
\end{lemma}
%

\subsection{Stochastic local stability implies global stability}
\label {sec:genglobloc} 

  Now we proceed to more elaborate situations, where local stability is provided by the noise perturbations of control, which was proved by the application of the Kolmogorov's Law of Large Numbers,  Lemma~\ref{lem:Kolm}.  
We follow the ideas of \cite{BR2019,Medv} and 
also \cite{BR0, BR1, BKR2020}. 
  
The proof of the main result, Theorem~\ref {thm:locglgen}, consists of two main  steps. On the first step, we show that, when $(\alpha, \ell)$ are chosen appropriately, a solution to \eqref{eq:PBCintr}  changes sides of $K$  at each step. Then we prove modified local stability:  for each $\nu\in (0, 1)$ we find a $\delta_0=\delta_0(\nu)>0$ s.t. as soon as  $x_t\in (K-\delta_0, K+\delta_0)$, $t\in \mathbb N$  is arbitrary, we get $\lim\limits_{n\to \infty} x_n=K$ on  the set which probability is not less than $1-\nu$. 
 The second step uses Borel-Cantelli Lemma~\ref{lem:topor} and is similar to the one in the proof of
Theorem~\ref{thm:locgldet}.
 
 We assume that \eqref {cond:loc2} holds with some $\theta>0$. In this section we concentrate on $f$ that  changes the side of $K$ at  each consecutive step, i.e.  
 $(f(x)-K)(K-x)>0$ in $(K-\theta, K+\theta)$,  which implies that the solution of equation \eqref{eq:originalintr} alternates its position relative to $K$ at each step. To guarantee this,  we assume that,  for some $a_1, a_2>0$,
 \begin{equation}
  \label{cond:flhdif}
   f(x)>a_1(K-x)+K, \quad x\in (K-\theta, K), \quad  f(x)<a_2(K-x)+K, \quad x\in (K, K+\theta).   
   \end{equation}
   We use constants $a_1$ and $a_2$ to impose assumptions on a control parameter  to guarantee that  the solution of equation \eqref{eq:PBCGstoch} also changes its position relative to $K$ at each step.   
   
Instead of assuming $\alpha>\Psi(\tilde L^-, \tilde L^+)$, where $\Psi$ was defined in \eqref{def:Psi}, as in the case of deterministic local stability, see Lemma~\ref{lem:locdet}, we introduce the condition
 \begin{equation}
 \label{cond:lambda0}
 \mathbb E\ln \left|{\mathcal L}^-(\alpha+\ell \xi) {\mathcal L}^+(\alpha+\ell \xi)\right| =
 -\lambda_0<0.
 \end{equation} 
 Here $\lambda_0>0$  is a positive number, and $ {\mathcal L}^{\pm}(\beta)=(1-\beta) \tilde L^{\pm}-\beta$.

Lemma~\ref{lem:locstoch} below states that conditions \eqref {cond:loc2}, \eqref{cond:flhdif},  
\eqref{cond:lambda0} and 
\begin{equation}
 \label{cond:sides}
 \alpha+\ell\le \min\left\{\frac{a_1}{a_1+1}, \, \frac{a_2}{a_2+1}\right\} 
 \end{equation}
guarantee local stochastic stability with any a priori given probability $1-\gamma$, $\gamma\in (0,1)$.  However, the smaller $\gamma$ is, the smaller $\delta_0$ in the local stability interval $(K-\delta_0, K+\delta_0)$ is required.

\begin{lemma}
\label{lem:locstoch}
Let  Assumptions  \ref{as:noise}, \ref{as:Lglob1} and conditions  \eqref {cond:loc2}, \eqref{cond:flhdif} hold,
and $(\alpha, \ell)$ satisfy  \eqref{cond:lambda0} and \eqref{cond:sides}. Then, for each $\gamma \in (0, 1)$  and $\delta_0>0$, there exists $\Omega_\gamma\subseteq \Omega$, $\mathbb P(\Omega_\gamma)\ge 1-\gamma$, s.t. for any solution $x$  to \eqref{eq:PBCintr} with $x_0\in (K-\delta_0, K+\delta_0)$, we have $\lim\limits_{n\to \infty} x_n=K$  on $\Omega_\gamma$.
\end{lemma}

Since Lemma~\ref{lem:locstoch} is a partial case of Step (ii) in the proof of the Theorem~\ref {thm:locglgen}, its proof is omitted. 

\begin{theorem}
 \label{thm:locglgen}
 Let  Assumptions  \ref{as:noise},\ref{as:Lglob1} and conditions  \eqref {cond:loc2}, \eqref{cond:flhdif} hold, 
 and  $\underline \beta_*$ be  defined in \eqref{def:beta*Gull}.
   Then $\displaystyle \lim_{n\to \infty} x_n=K$  a.s., for any solution $x$  to \eqref{eq:PBCintr} with $x_0>0$ and  $(\alpha, \ell)$ satisfying  \eqref{cond:lambda0}, \eqref{cond:sides} and 
 \begin{equation}
 \label{cond:glob}
 \alpha\in (0, \, \underline\beta_*), \,\, \ell \in \left( \underline\beta_*-\alpha, \min\{\alpha, 1-\alpha  \}  \right).
  \end{equation}
   \end{theorem}

 \begin{remark}
 \label{rem:stochbetter}
 It is straightforward to check  that \eqref{cond:lambda0} is satisfied if $\alpha>\underline \beta_0=\Psi(\tilde L^-, \tilde L^+)$, $\ell=0$. For Bernoulli distributed noises $\xi_n$, it can be simply  shown that there are $(\alpha, \ell)$, $\alpha<\underline \beta_0$  and $\ell>0$ s.t. \eqref{cond:lambda0} holds. 
Indeed, in the case of Bernoulli distributed $\xi$, 
  \begin{equation}
  \begin{split}
  \label{def:mathcalV}
  &\mathbb E\ln \left[{\mathcal L}^-(\alpha+\ell \xi) {\mathcal L}^+(\alpha+\ell \xi)\right]=\frac 14\ln\left(\mathcal V(\alpha, \ell)\right),\\
  &\mathcal V(\alpha, \ell):=\left[(L^--\alpha(L^-+1))^2-\ell^2(L^-+1)^2  \right]\left[(L^+-\alpha(L^++1))^2-\ell^2(L^++1)^2  \right].
\end{split}
\end{equation}
By Lemma~\ref{lem:mathcalL}~(vi), ${\mathcal L}^-(\underline \beta_0) {\mathcal L}^+(\underline \beta_0)=1$, so
$\mathcal V(\underline \beta_0, 0)=1$. Then there exists $\ell>0$ s.t. each bracket on the second line of  \eqref {def:mathcalV} becomes smaller but remains positive, so  $\mathcal V(\underline \beta_0, \ell)<1$. Now we can choose $\beta_1\in (0, \underline \beta_0)$ s.t. $\mathcal V(\beta_1, \ell)<1$. Thus by introducing noise into control, we can allow smaller average control values, see more details for continuously differentiable unimodal $f$ 
in Section~\ref{subsec:Zinger}. 
\end{remark}
  
 \subsection{Unimodal continuously differential $f$: when local implies global stability}  
 \label{subsec:Zinger} 
To ensure equivalence of local and global deterministic stability, we impose additional restrictions on $f$.
\begin{assumption}
\label{as:Singer}
Let $f$ satisfy Assumption~\ref{as:Lglob1}, be unimodal three times differentiable with a unique critical point $c\in (0,K)$ (maximum), $f'(0)>1$, $f''(x)<0$ for all $x \in (0,c)$ and a negative Schwarzian derivative \eqref{eq:Schwarzian} everywhere but at $c$.  
\end{assumption}
If $f$ satisfies Assumption~\ref{as:Singer}, local stability of $K$ implies the global one.
Moreover, this is also true for PBC of $f$ \cite{FL2010}.
Under Assumption~\ref{as:Singer}, the value  of the parameter $\displaystyle \alpha > \alpha_0 := \frac{-f'(K)-1}{-f'(K)+1}$ provides local, as well as global stability of \eqref{eq:PBCintrdet}, see \cite {FL2010}, we can choose $\underline \beta_*=\alpha_0$. 
If $\alpha < \alpha_0$, we get $\displaystyle \frac{d}{dx} G(\alpha,K)<-1$, which means, due to smoothness, that in some neighbourhood of $K$, we get $(x-K)(G(\alpha,x)-K)<0$ and $|G(\alpha,x)-K| > |x-K|$, thus  $K$ is repelling solutions in a certain neighbourhood. 
 
In this section we show that, for each $f$ satisfying Assumption~\ref{as:Singer}, in the cases of Bernoulli and uniformly distributed noises $\xi$, we can decrease the mean value  to $\alpha<\alpha_0$, so that the stochastic control $\alpha+\ell \xi_n$  with a specially chosen $\ell< \min\{\alpha, 1-\alpha\}$, provides global stability of the solution $x$ to \eqref{eq:PBCintr} with probability one. 

By \cite{FL2010}, we can take $\underline \beta_*=\alpha_0$. We only consider $\displaystyle\alpha<\frac{-f'(K)-1}{-f'(K)+1}$, which implies $\mathcal L(\alpha)>1$ (otherwise, we get stability without noise).  Set 
\begin{equation}
\label{def:L0a0}
L_0:=-f'(K), \quad  \alpha_0 := \frac{-f'(K)-1}{-f'(K)+1}.
\end{equation}
\begin{remark}
\label{rem:sharpcond}
Further, we use the expression 
$\mathbb E\ln \left[(1-\alpha-\ell \xi) L_0-\alpha-\ell \xi\right]$ which should be negative for local stability.
If $\ell=0$, we get  $0<(1-\alpha)L_0-\alpha < 1$, or $\alpha \in (\alpha_0,1)$, 
where $\alpha_0$ is from \eqref{def:L0a0}, which is the well-known sharp stability condition in the deterministic case \cite{FL2010}.
\end{remark}

\begin{theorem}
\label{thm:Singer}	
 Let  Assumptions  \ref{as:noise},  \ref{as:Singer} hold, $L_0$ and  $\displaystyle \alpha_0$ be defined as in \eqref{def:L0a0}. Assume that,  for some $\alpha \in (0,\alpha_0)$, $\ell \in \left( \alpha_0-\alpha, \min\{ \alpha, 1-\alpha \} \right)$, $\lambda_1>0$, either 
\begin{equation}
\label{cond:lambda11}
\alpha+\ell< \frac {L_0}{L_0+1}, \quad
\mathbb E\ln \left[(1-\alpha-\ell \xi)L_0 -\alpha-\ell \xi\right]=-\lambda_1<0
\end{equation}
or 
\begin{equation}
\label{cond:lambda1}
\mathbb E\max\left\{  \ln \left[(1-\alpha-\ell \xi)L_0\right], \,   \ln[\alpha+\ell \xi]\right\}=-\lambda_1<0
\end{equation}
holds. Then $\displaystyle \lim_{n\to \infty} x_n=K$  a.s., for any solution $x$  to \eqref{eq:PBCintr} with $x_0>0$.
\end{theorem}
      
The next theorem demonstrates that  when $f$ satisfies Assumption~\ref{as:Singer}  and  $\xi$  has a symmetric  distribution, we can decrease the average value of the control, compared to the minimal deterministic one $\alpha_0$.
\begin{theorem}
\label{thm:contdiscrdistr}
Let  $f$ satisfy Assumption \ref{as:Singer},  $L_0$ and $\alpha_0$  be defined as in \eqref{def:L0a0}, 
Assumption~\ref{as:noise} hold, random variables $\xi_n$ have a symmetric (around zero) distribution, and be either continuous or discrete with a countable number of states. Then there exist $\alpha \in (0,\alpha_0)$ and  $\ell \in \left( \alpha_0-\alpha, \min\{ \alpha, 1-\alpha \} \right)$ s.t.  $\alpha_0<\alpha+\ell<\frac{L_0}{L_0+1}$ and \eqref{cond:lambda11} holds.
\end{theorem}

\begin{remark}
\label{rem:bernuncont}
Note that the values of $\alpha<\alpha_0$ and $\ell$ established in Theorem~\ref{thm:contdiscrdistr}  are not supposed to be optimal (for example, the minimal $\alpha$), we just show that they exist for either continuous or discrete distribution of $\xi$. For each  particular distribution  we can find smaller values of $\alpha$ by calculating $\mathbb E \ln [\mathcal L_0(\alpha+\ell\xi)]=\mathbb E \ln [L_0-(\alpha+\ell\xi)(L_0+1)]$.
 
When $\xi$ has a Bernoulli distribution, we get  $\mathbb E \ln [\mathcal L_0(\alpha+\ell\xi)]=\frac 12\ln \left [[L_0-\alpha(L_0+1)]^2-\ell^2(L_0+1)^2 \right],$ so \eqref{cond:lambda11} holds when
\begin{equation}
\label{cond:bern}
 \left[\frac {L_0}{L_0+1} -\alpha\right]^2-\frac1{(L_0+1)^2}   <\ell^2<\min\left\{\left[\frac {L_0}{L_0+1} -\alpha\right]^2, \alpha\right\}.
 \end{equation}
For example, for $L_0>2$, the values $\displaystyle \alpha=\frac {L_0-1-\frac{L_0-2}2}{L_0+1}$,  $\displaystyle \ell=\frac {\frac{L_0-2}2+ 0.5}{L_0+1}$ satisfy $\alpha \in (0,\alpha_0)$ and \eqref{cond:bern}, a similar example can be found for $L_0\in (1,2]$.

When $\xi$ has the uniform continuous on $[-1,1]$ distribution,  
we obtain $\mathbb E \ln [\mathcal L_0(\alpha+\ell\xi)]=\frac 12\int_{-1}^1\ln  [\mathcal L_0(\alpha+\ell u)]du$, so \eqref{cond:lambda11} holds when
\begin{equation}
\label{cond:cont}
\ln [\mathcal L_0(\alpha)]<\frac1{6}\left[\frac{\ell(L_0+1)}{\mathcal L_0(\alpha)} \right]^2.
\end{equation}
Using estimation of the integral, example pairs $(\alpha,\ell)$ can be found with $\alpha \in (0,\alpha_0)$ for which
\eqref{cond:cont} is valid.

Example~\ref{ex:Rick} considers the Ricker model with a control perturbed by the  Bernoulli
or the continuous noise.  
The simulation results for the  Bernoulli perturbations illustrate that parameters computed using \eqref{cond:bern} are quite sharp. 
Similar calculations can be implemented for the uniform distribution illustrating sharpness of \eqref{cond:cont}.
\end{remark}

 
 \section{Determining control of the deterministic equation}
\label{subsec:parglobstab}

In this section we discuss situations when we are able to find the parameter $\beta_*$ which guarantees global stability of the solution to deterministic equation \eqref{eq:PBCintrdet}. In all results $\beta_*$ might be not optimal (the Ricker model demonstrates this, see Example \ref{ex:Rick}) even though it is much better than controls found in  \cite{BKR2016, BKR2020} based on the global constants.

We are going to use the method of envelope  functions suggested by Cull~\cite{Cull2}. 
The next result from \cite[Theorem 3, P. 996]{Cull2} is slightly adapted to our needs. 

\begin{lemma}
\label{lem:Cull2}
Let $\phi(x)$ be a monotone decreasing function which  is positive on $(d_1, d_2)$ and $\phi(\phi(x))=x$, $K \in (d_1,d_2)$. 
Assume that $f(x)$ is a continuous function s.t. $\phi(x)>f(x)$ for $x\in (d_1, K)$, 
$\phi(x)<f(x)$ for $x\in (K, d_2)$, $f(x)<x$ for $x>K$, $f(x)>x$ for $x\in (d_1, K)$, 
$f(x)>0$ on $(d_1, d_2)$. Then, for all $x\in (d_1, d_2)$, $\lim\limits_{k\to \infty}f^{(k)}(x)=K$.
\end{lemma}

When $\phi$ is such as in Lemma~\ref{lem:Cull2}, we say that $\phi$ envelopes $f$, see  \cite{Cull2}.

We start with the case when the function $f$ satisfies only Assumption \ref{as:Lglob2}.

\begin{proposition}
\label{prop:Globdet}
Let  Assumption~\ref{as:Lglob2}  hold, $1<L^+\le L^-$, $\underline \beta_0 $ be defined as in  \eqref{def:beta0}.  Then $\lim\limits_{n\to \infty}x_n=K$ for each solution $x$ to \eqref{eq:PBCintrdet} with $\alpha\in (\underline \beta_0, 1)$ and $x_0>0$.
\end{proposition}

Now we generalize Proposition~\ref{prop:Globdet} to the case when the interval $[x_m, f_m]$ is split into several subintervals and $f$ satisfies a Lipschitz condition with different constants on each of them. In many population dynamics models, such as Ricker's and logistic, using just the left-hand-side $L^-$ and the right-hand-side $L^+$ Lipschitz type constants does not give the best possible value for the control parameter. Models like Ricker's have the property  that to the left of $K$ the Lipschitz-type constants are much larger than the derivative at $K$, while to the right of $K$ they are quite small. In the following  we are going to use such situation and find a better low bound for the control than $\underline \beta_0$.  
First,  we construct  a piecewise function which  can be used as $\phi$ in Lemma~\ref{lem:Cull2}
for corresponding  function  $G(\alpha, x)$.  Consider  finite sequences $(a_i)_{i=0, 1, \dots, m}$, $(C^-_i)_{i=0,\dots, m-1}$, $(C^+_i)_{i=0, \dots, m-1}$, $(b_i)_{i=0, \dots, m-1}$,
\begin{equation}
\label {def:seq1}
\begin{split}
& 0<a_m<\dots a_2<a_1<a_0 := K,
\\
& C_i^-, C_i^+>0, \, \, C_i^-C_i^+=1,\quad
 b_{i+1}=-C_i^-(a_{i+1}-a_i)+b_i, \,\,  i=0, \dots, m-1, \,\, b_0 := K,
\end{split}
\end{equation}
 and a piecewise function $\phi:(0, \infty)\to (0, \infty)$, $\phi(K)=K$, is defined for $i=0,\dots m-1$ by
\begin{equation}
\label {def:phi}
\begin{split}
&\phi(x)=-C_0^-(x-K)+K, \quad a_1\le x\le K, \quad \phi(x)=-C_0^+(x-K)+K, \quad K\le x\le b_1,\\
&\phi(x)=-C_1^-(x-a_1)+b_1, \quad a_2\le x < a_1, \quad \phi(x)=-C_1^+(x-b_1)+a_1, \quad b_1 < x\le b_2,\\
&\phi(x)=-C_i^-(x-a_i)+b_i, \, a_{i+1}\le x < a_i, \quad \phi(x)=-C_i^+(x-b_i)+a_i, \,  b_i <x\le b_{i+1}.
\end{split}
\end{equation}
For each $i=0, \dots, m-1$, we have $b_{i+1}=\phi(a_{i+1})$. Also, when $a_{i+1}\le x\le a_i$,
\[
\phi(\phi(x))=-C_i^+(C_i^-(x-a_i)+b_i-b_i)+a_i=x,
\]
and,  when  $b_{i}\le x\le b_{i+1}$,
\[
\phi(\phi(x))=-C_i^-(C_i^+(x-b_i)+a_i-a_i)+b_i=x.
\]
So $\phi$ is a monotone decreasing function which is positive on $(0, b_m)$, since $\phi(b_m)=a_m$ and $\phi(\phi(x))=x$.
Now we proceed to equation \eqref{eq:PBCintrdet}. Let \eqref{def:seq1} hold, $L_i^-, L_i^+>0, $ $i=0, \dots, m-1$,  $L_0^-L_0^+>1$,
\begin{equation}
\label{ineq:1group}
f(x)-K\le L_0^-(K-x),\,\,  a_{1}\le x\le K,\quad K-f(x)\le L_0^+(K-x), \,\, x>K,
\end{equation}
\begin{equation}
\label{def:alpha00}
\alpha_0 =\Psi \left(L_0^-, L_0^+ \right),
\end{equation}
and, for $ i=1, \dots, m, $ $\mathcal L^{\pm}_i(\alpha):=L_i^{\pm }-\alpha(L_i^{\pm }+1)$, 
\begin{equation}
\label {cond:fL}
\begin{split}
 &b_{i}(\alpha_0)=\mathcal L_{i-1}^-(\alpha_0)(a_{i}-a_{i-1})+b_{i-1}(\alpha_0),\quad b_0:=K,\\
&f(x)-f(a_i)\le L_i^-(a_i-x),\, a_{i+1}\le x\le a_i,   \quad f(b_i(\alpha_0))-f(x)\le L_i^+(x-b_i(\alpha_0)),\,\,  b_{i}(\alpha_0)\le x. 
\end{split}
\end{equation}

\begin{remark}
\label{rem:rem:condmultintr}
We assume $L_0^-L_0^+>1$ since otherwise the equilibrium $K$ is locally stable. We also assume without loss of generality that 
 $\mathcal L_{i-1}^-(\alpha_0)>0$ for all $i=1, \dots, m$. Indeed, if, for some $i$, we have $\mathcal L_{i-1}^-(\alpha_0)<0$, it means that $b_{i}(\alpha_0)<b_{i-1}(\alpha_0)$, so we exclude  $a_i$ from the sequence in~\eqref{def:seq1}.
\end{remark}

\begin{proposition}
\label{prop:envG}
Let Assumption~\ref{as:Lglob2}  and conditions \eqref{def:seq1}, \eqref{ineq:1group}-\eqref{cond:fL} hold and also 
\begin{equation}
\label{cond:alpha0gl}
\alpha_0 \ge  \max_{i=1,\dots, m}\left\{\frac{L_i^-L_i^+-1}{(L_i^-+1)(L_i^++1)} \right\}.
\end{equation}
If $x_n$  is a solution to \eqref{eq:PBCintrdet} with any $\alpha >\alpha_0$ and $x_0>0$  then  $\lim\limits_{n\to \infty}x_n=K$.
\end{proposition}
 Proposition \ref{prop:envG} is illustrated in Example~\ref{ex:local_is_global}~(a).

\begin{remark}
\label{cor:glob}
By Lemma~\ref{lem:locdet}, conditions \eqref{ineq:1group}-\eqref{def:alpha00} provide stability of the solution to \eqref{eq:PBCintrdet} with $\alpha >\alpha_0$ on the interval $\bigl(a_1, \mathcal L_{0}^-(\alpha_0)(a_{1}-K)+K\bigr)$. Therefore, Proposition \ref{prop:envG} proves that when conditions \eqref {cond:fL}-\eqref{cond:alpha0gl} hold,  the local stability implies the global one.
\end{remark}

We also can deal  with the case when there is $i_0<m$ s.t.  $\alpha_0<\Psi \left(L_{i_0}^-, L_{i_0}^+ \right)$ and find a bigger parameter $\bar \alpha$ which guarantees global stability of the solution to \eqref{eq:PBCintrdet}.  It is discussed in Remark \ref {rem:glob} in the Appendix and illustrated in Example~\ref{ex:local_is_global}~(b).

When $f$ is continuously differentiable, we can get an explicit result computing  $\alpha_0$ which ensures global stability. 

\begin{assumption}
\label{as:3}
Let  $f$ satisfy Assumption \ref{as:Lglob1}, be continuously differentiable with  $f'(x)<1$ for  $ x\in [x_{\max}, f_m]$, and for $\alpha_0$ defined in \eqref{def:L0a0},  
\begin{equation}
\label{cond:alpha0}
\alpha_0>\Psi \biggl(f'(G(\alpha_0, x)), f'(x) \biggr) \quad \mbox{for each} \,\, x\in [x_{\max}, K).
\end{equation}
\end{assumption}

Note that since Assumption \ref{as:Lglob1} holds, there is no local stability at $K$, and therefore $f'(K)<-1$,
so $\alpha_0$ is well defined. Also, any $\alpha\in (\alpha_{0}, 1)$ leads to local stability of $K$ for \eqref{eq:PBCintrdet}.

\begin{proposition}
\label{prop:3}
Let Assumption \ref{as:3} hold, and $\alpha\in (\alpha_{0}, 1)$. Then $\lim\limits_{n\to \infty}x_n=K$, where $x$ is a solution to \eqref{eq:PBCintrdet} with  $x_0\in (0, \infty)$. 
\end{proposition}

\begin{figure}[ht]
\centering
\vspace{-60mm}
\includegraphics[height=.5\textheight]{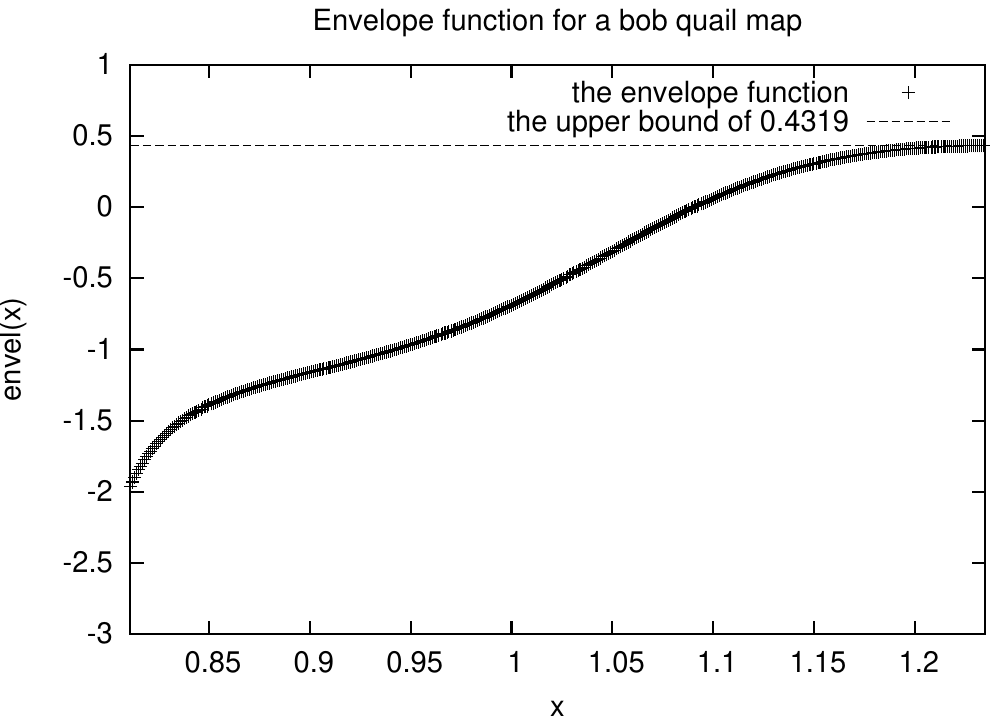}
\caption{The envelope function $\displaystyle {\rm envel}(x) := \Psi\left(g'(x), g'(G( \tilde \alpha, x)) \right)$ 
for $x\in [x_{\max}, K] \approx [0.81,1.2347]$ compared to the bound for the control $\alpha$  which is $\tilde \alpha \approx 0.4319$.
We observe that $ \tilde \alpha\ge \max_{x\in [x_{\max}, K]}{\rm envel}(x) $ 
leading to both local and global stability for $\alpha \in (\tilde{\alpha},1)$.}
\label{figure12}
\end{figure}

To illustrate Proposition~\ref{prop:3}, we consider the bobwhite quail map \cite{mb}
 \begin{equation}
 \label{def:Bobq}
 g(x)=x\left(0.55+\frac{3.45}{1+x^9} \right), \quad x>0,
  \end{equation}
where $K \approx 1.2347$, $f'(K) \approx  2.521$, $\tilde \alpha  \approx \frac{1.521}{3.521} \approx 0.4319$, $ x_{\max} \approx 0.811$, $\Psi(a,b)=\frac {ab-1}{(a+1)(b+1)}$. 
Fig.~\ref{figure12} shows that $ \tilde \alpha\ge \max_{x\in [x_{\max}, K]}\Psi\left(g'(x), g'(G( \tilde \alpha, x)) \right) = \max_{x\in [x_{\max}, K]} {\rm envel}(x) $.

\section{Examples and simulations}
\label{sec:exsim}

 \begin{example}  
 \label{ex:Rick} 
Consider Ricker's function $ f(x)=xe^{r(1-x)}$, $x\ge 0$, which satisfies
$\underline \beta_*=\alpha_0=\frac {-f'(1)-1}{-f'(1)+1}$.

(a) Let $r=3.5$, then  $-L_0=f'(1)=1-r=-2.5$, $\alpha_0=\underline \beta_*=3/7 \approx 0.4285$.  In the case of Bernoulli distributed $\xi$
we apply  formula \eqref{cond:bern} from Remark \ref{rem:bernuncont}  and get that $\ell$ should satisfy $\sqrt{(0.71-\alpha)^2-0.0811}<\ell<\min\{(0.71-\alpha), \alpha \}$, see the domain in Fig.~\ref{figure1}. Taking $\alpha=0.368$ we should have $\ell \in (0.1877, 0.342)$.  
This case is illustrated by the bifurcation diagram in Fig.~\ref{figure2}, left. The runs for $\ell =0.2$ and $\alpha = 0.3,0.36,0.37$ in Fig.~\ref{figure3} also confirm this.
 
For uniformly distributed $\xi$ we apply inequality  \eqref{cond:cont} from Remark \ref{rem:bernuncont}, which is satisfied  when $\alpha=0.405$ and $\ell=0.2$, which coincides with what we observe on  the bifurcation diagram in Fig.~\ref{figure2}, right. 
\begin{figure}[ht]
\centering
\vspace{-60mm}
\includegraphics[height=.5\textheight]{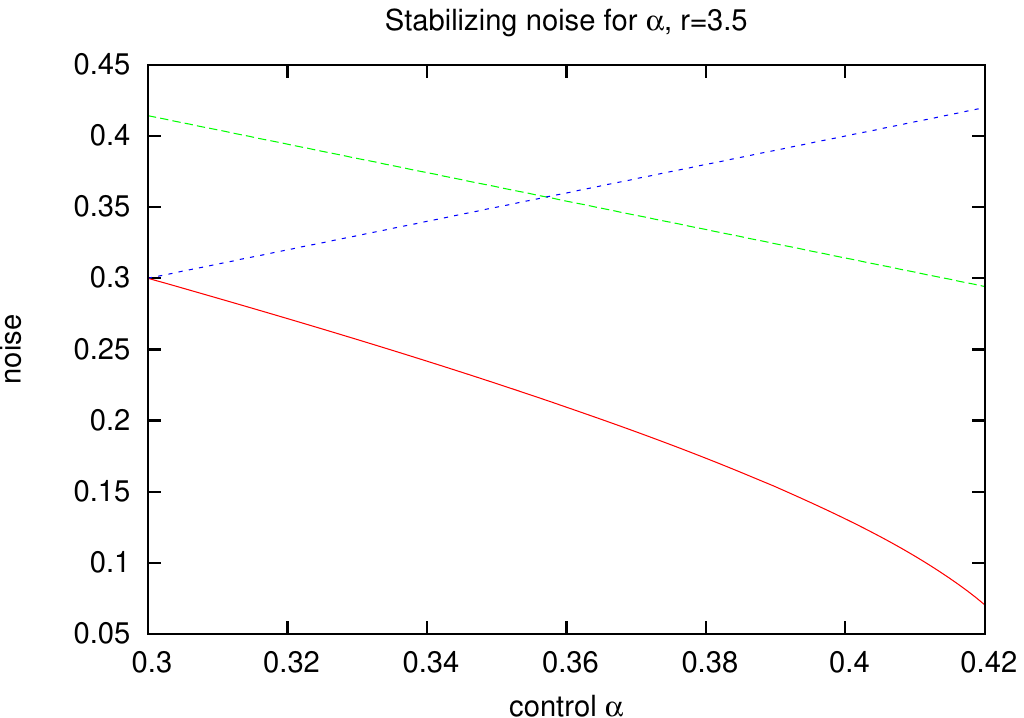}
\caption{The domain of the Bernoulli noise $\ell$ leading to stability for the Ricker model with $r=3.5$ and PBC 
with a given $\alpha$ is above the red solid line. 
In addition, $\ell < \alpha$ limits the allowed values below the quadrant bisect (the dashed blue line).
}
\label{figure1}
\end{figure}

\begin{figure}[ht]
\centering
\vspace{-45mm}
\includegraphics[height=.4\textheight]{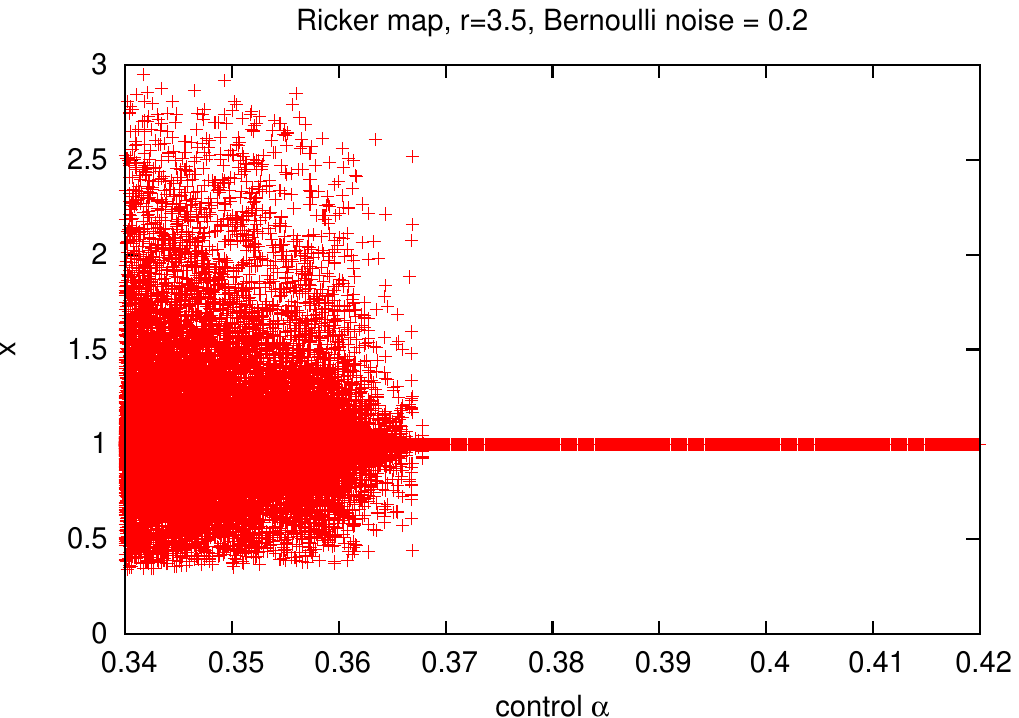}
~~
\includegraphics[height=.4\textheight]{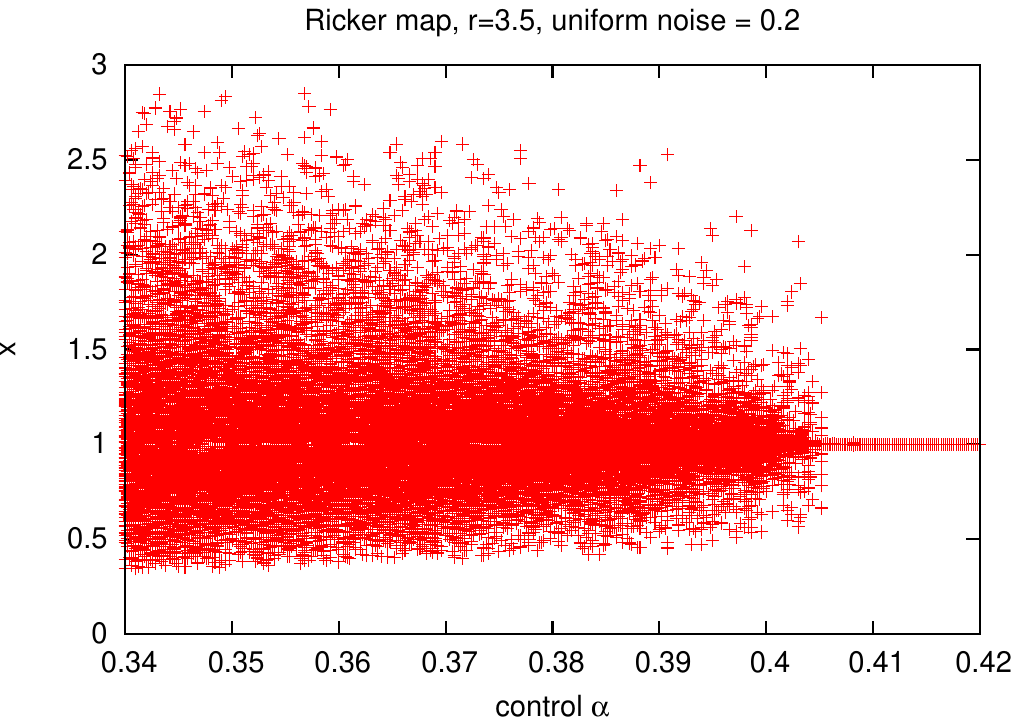}
\caption{Bifulcation diagrams for the Ricker map stabilized by PBC with $r=3.5$ and $\alpha$ perturbed by (left) the Bernoulli 
noise taking values $\pm 0.2$; (right) uniformly distributed in $[-0.2,0.2]$ noise.} 
\label{figure2}
\end{figure}

\begin{figure}[ht]
\centering
\vspace{-40mm}
\includegraphics[height=.3\textheight]{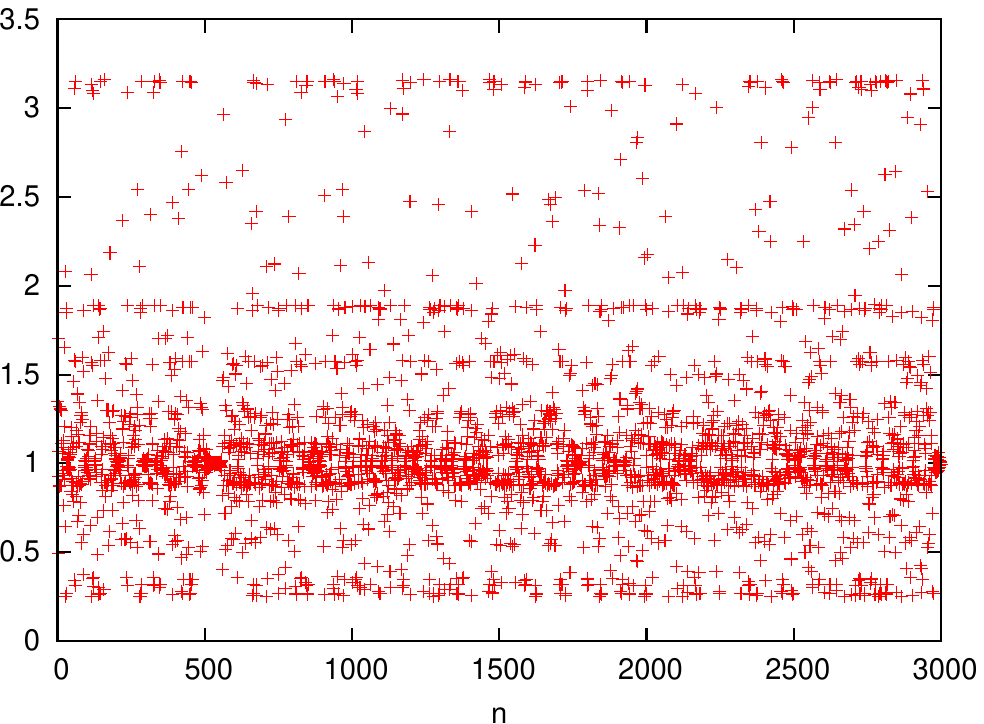}~~~
\includegraphics[height=.3\textheight]{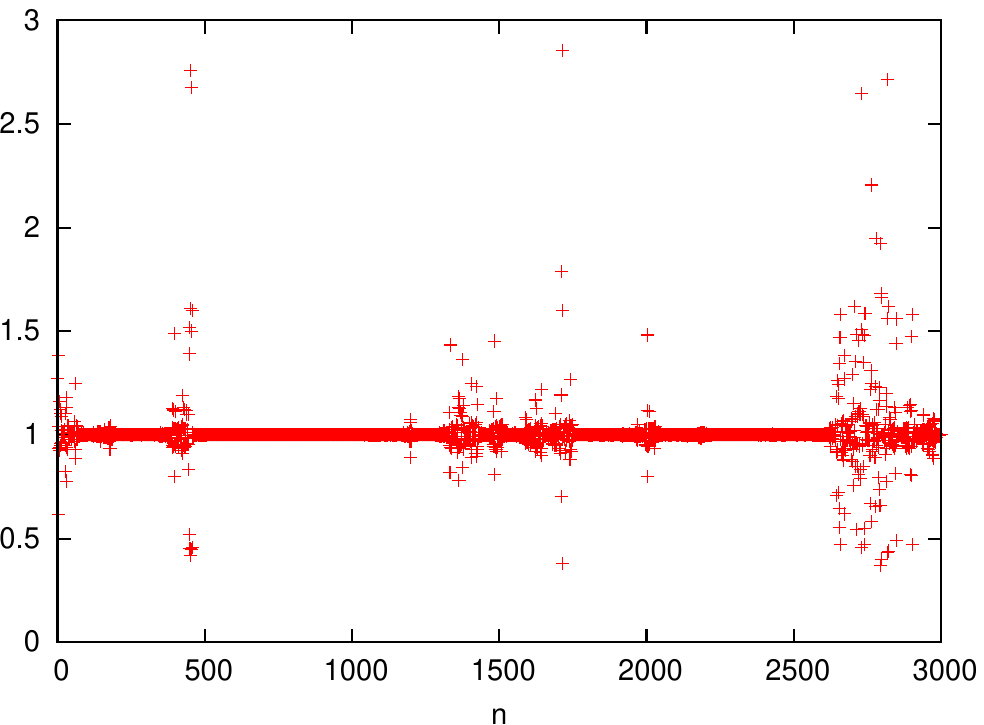}~~~
\includegraphics[height=.3\textheight]{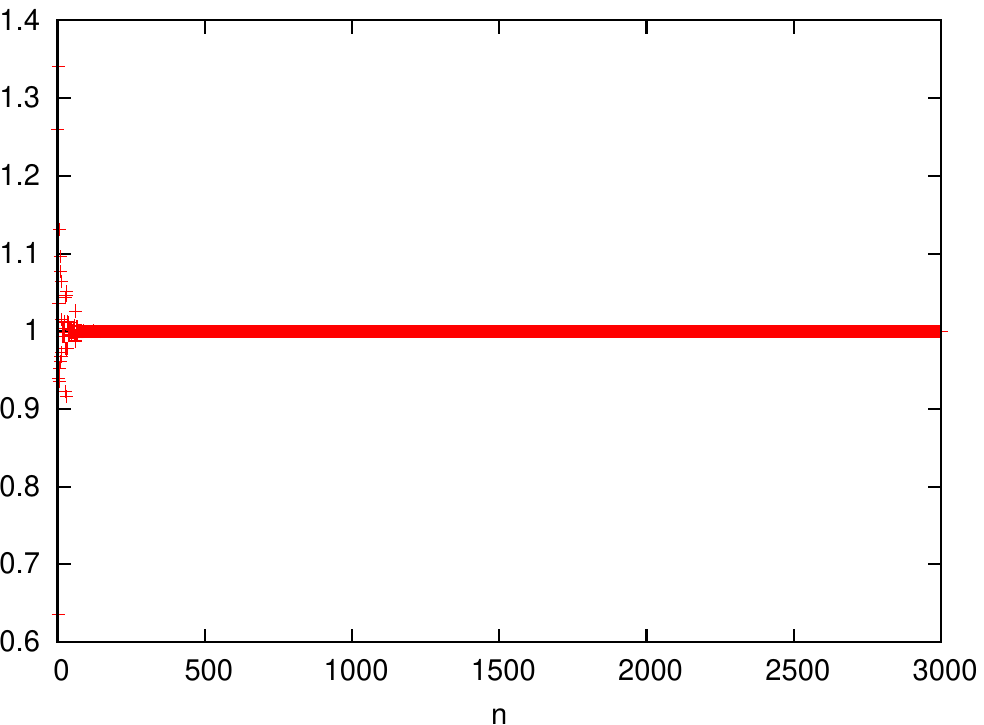}
\caption{Runs of the controlled Ricker model with $\alpha=3.5$, the Bernoulli
noise taking values $\pm 0.2$ and (left) $\alpha=0.3$, (middle) $\alpha=0.36$, (right) $\alpha=0.37$.}  
\label{figure3}
\end{figure}

(b)  Let  now $r=3$,  then $-L_0=f'(1)=1-r=-2$, $\alpha_0=1/3=0.33(3)$. For Bernoulli distributed $\xi$, 
the domain for $(\alpha,\ell)$ based on condition \eqref{cond:bern} from Remark~\ref{rem:bernuncont}, is described in Fig.~\ref{figure4}. We conclude that for $\alpha>0.283$, $\ell=0.2$, the equilibrium $K$ is globally stable. 
The bifurcation diagram for $\ell = 0.2$ on Fig.~\ref{figure5} confirms it.

\begin{figure}[ht]
\centering
\vspace{-60mm}
\includegraphics[height=.5\textheight]{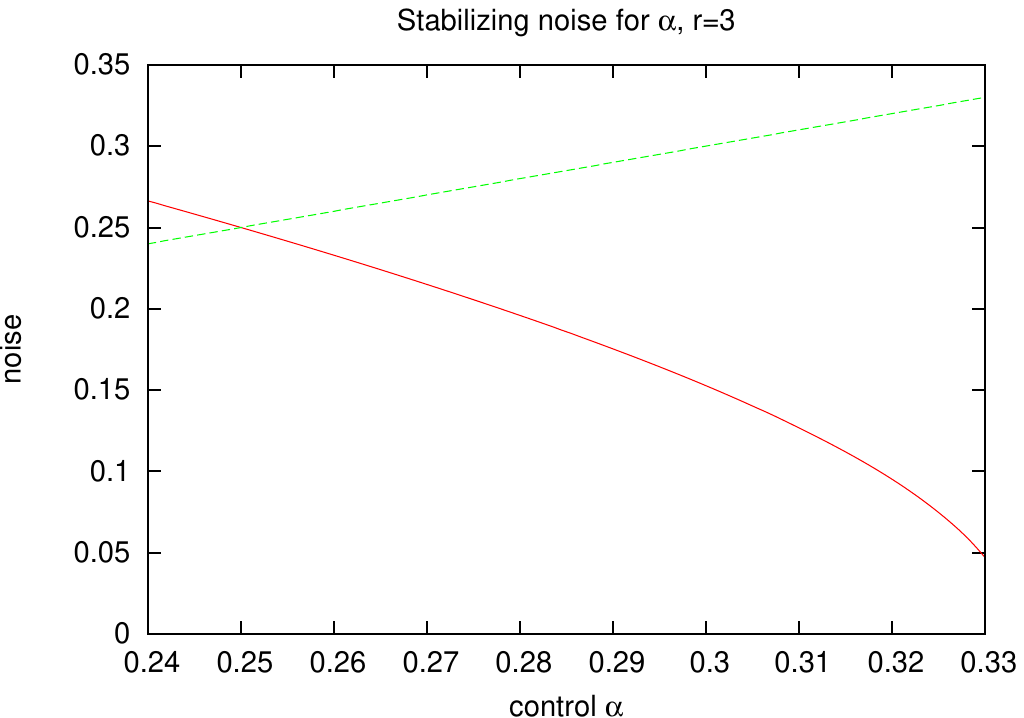}
\caption{The domain of the Bernoulli noise $\ell$ leading to stability for the Ricker model with $r=3.0$ 
and PBC with a given $\alpha$ is above the red solid line and below the green dashed line. 
The upper bound is above the bisect $\ell = \alpha$ marking that the noise cannot exceed the control (the dashed green line) and thus is not relevant in this case.
}
\label{figure4}
\end{figure}

\begin{figure}[ht]
\centering
\vspace{-60mm}
\includegraphics[height=.5\textheight]{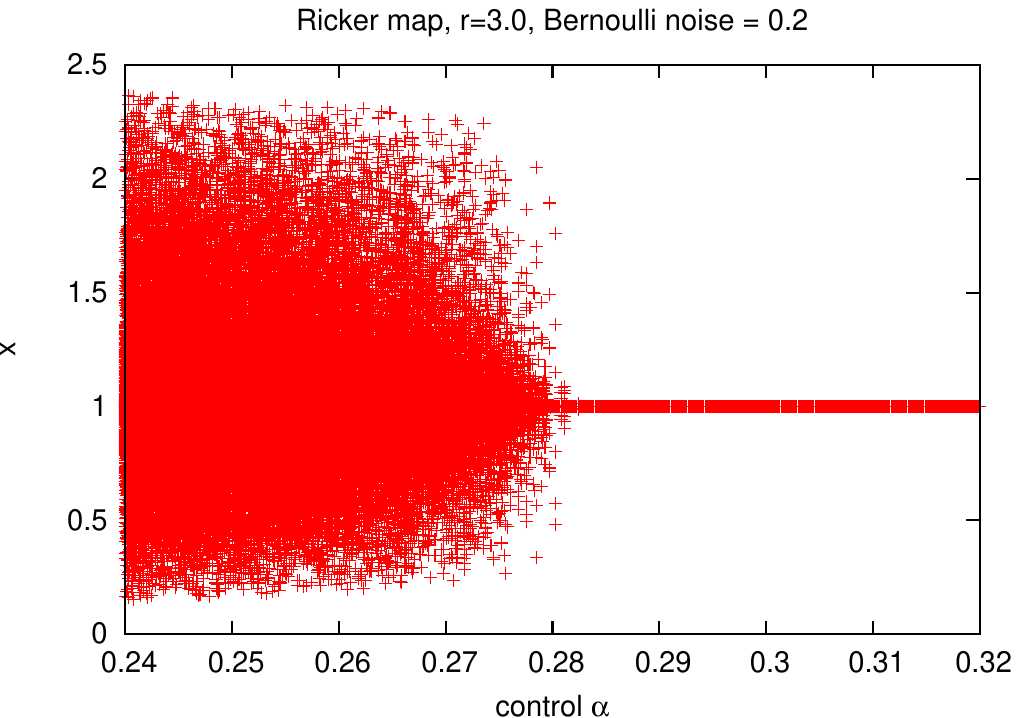}
\caption{Bifulcation diagrams for the Ricker map with $r=3.0$ stabilized by PBC with $\alpha$ perturbed by the Bernoulli 
noise taking values $\pm 0.2$} 
\label{figure5}
\end{figure}
Simulations illustrate sharpness of our theoretical computations: the final bifurcation leading to stability in Fig.~\ref{figure5} is quite close to theoretically computed $\alpha=0.283$. 
Applying Figs.~\ref{figure2} and~\ref{figure5}, or calculating directly by formulas \eqref{cond:bern}, \eqref{cond:cont}, we can find smaller average stabilizing value $\alpha$ of the stochastic control, but it might involve larger noise intensity $\ell$.
\end{example}

\begin{example}
\label{ex:local_is_global}
 (a)  Consider a piece-wise linear function with a unique positive fixed point $K=32$
   \begin{equation}
\label{ex:glob}
f(x):=\left\{\begin{array}{ll}
   51x/28 , & x\in [0, 28),\\
   -6x+219, & x\in [28, 29)=[a_3, a_2), \\  
   -5x+190, & x\in [29, 31)=[a_2, a_1), \\
   -3x+128, & x\in [31, 32)=[a_1, K), \\ 
    -2x+96, & x\in [32, 33) =[K, b_1), \\ 
  -1.4x+76.2, & x\in [33, 38)=[b_1, b_2), \\
    -1.2x+68.6, & x\in [38, 50)=[b_2, b_3), \\ 
      8.6, & x\ge 50.
     \end{array}\right.
      \end{equation}
  For function $f$ defined by \eqref{ex:glob}, Assumptions \ref{as:Lglob2} and conditions  \eqref{def:seq1},\eqref{ineq:1group}-\eqref{cond:fL} hold with   $K=32$, $x_{\max}=28$, $f_m=51$,  
	$a_1=31$, $a_2= 29 $, $a_3= 28=x_{\max}$,  $b_1=33$, $b_2= 38 $, $b_3= 39$, 
  \begin{equation*}
  \begin{split}
 & L^-_0=3, \,\, L_0^+=2, \,\,  \alpha_0=\Psi(2,3)\approx 0.416(6),\\
 &  L^-_1=5,\,\, L_1^+=1.4, \,\,\Psi(5, 1.4)\approx 0.416(6),\quad  L^-_2=6,\,\, L_2^+=1.2, \,\,\Psi(6, 1.2)\approx 0.4025.
 \end{split}
 \end{equation*}
 Since $\alpha_0=\Psi(2,3)\ge \max\{\Psi(5, 1.4), \Psi(6, 1.2)\}$, condition \eqref{cond:alpha0gl} holds and we can apply Proposition \ref{prop:envG} and conclude that $\lim\limits_{n\to \infty}x_n=K$ for any $\alpha >\alpha_0$ and $x_0>0$. If, for calculating global stability control,  we use maximum of left and right Lipschitz type constants with respect to the equilibrium $K$, $L^-=\frac{51-32}{32-28}=4.75$ and $L^+=2$ we arrive at 
 $\Psi(4.75, 2)=0.49275>0.4166=\alpha_0$, which shows the advantage of Proposition \ref{prop:envG} comparing  to Proposition \ref{prop:Globdet}.   

 Now we perturb the control with the noise $\ell \xi$, where $\xi$  has a Bernoulli distribution, and apply Theorem  \ref{thm:locglgen}, which allows to  decrease the control average $\alpha$. Based on the above, we need only to show that  $ \alpha+\ell>\alpha_0,$ and check condition \eqref{cond:lambda0} with
  ${\mathcal L}^-(\alpha+\ell \xi)=L^-_0-(\alpha+\ell \xi)(L^-_0+1)=3-4(\alpha+\ell \xi)$, 
  ${\mathcal L}^+(\alpha+\ell \xi)=L^-_0-(\alpha+\ell \xi)(L^-_0+1)=2-3(\alpha+\ell \xi)$.   Condition \eqref{cond:lambda0}    holds if 
  \begin{equation*} 
  \begin{split}
     &\mathcal V(\alpha, \ell):=\left[(L^-_1-\alpha(L^-_1+1))^2-\ell^2(L^-_1+1)^2  \right]\left[(L^+_1-\alpha(L^+_1+1))^2-\ell^2(L^+_1+1)^2  \right]\\
  &=
  \left[(3-4\alpha)^2-16\ell^2  \right]\left[(2-3\alpha)^2-9\ell^2\right] <1. 
  \end{split}  
   \end{equation*}   
   By direct calculations we show  that, when $\alpha=0.36, \ell=0.2$ we have  $ \alpha+\ell=0.56>0.417>\alpha_0,$ and    $
   {\mathcal V}(0.36, 0.2)\approx  0.8724<1$.

\begin{figure}[ht]
\centering
\vspace{-20mm}
\includegraphics[height=.23\textheight]{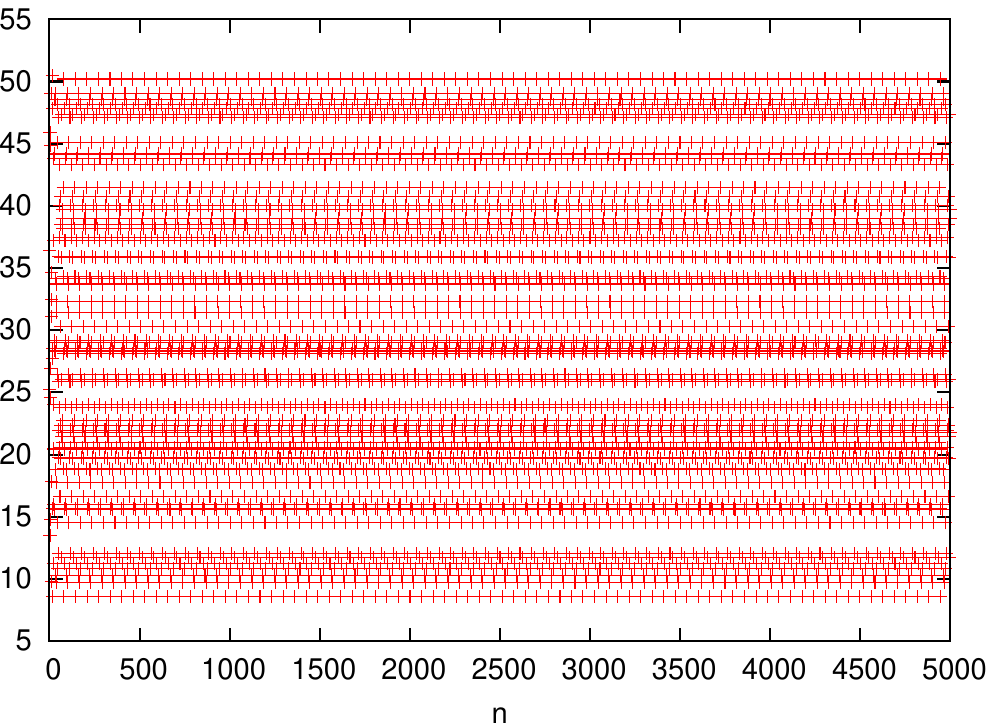}
\includegraphics[height=.23\textheight]{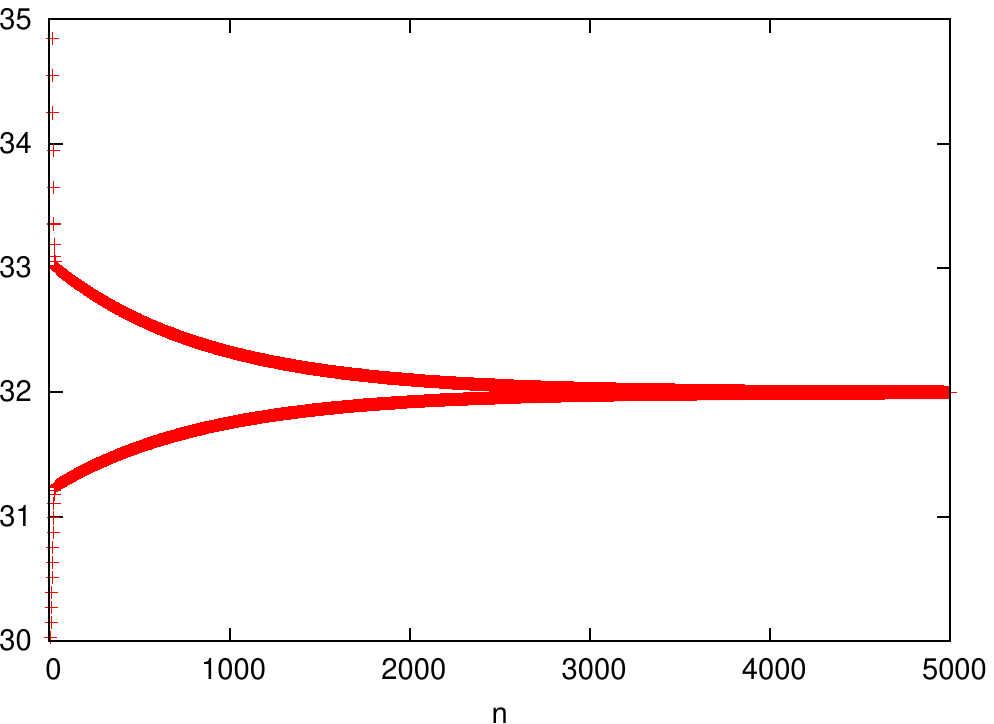}
\includegraphics[height=.23\textheight]{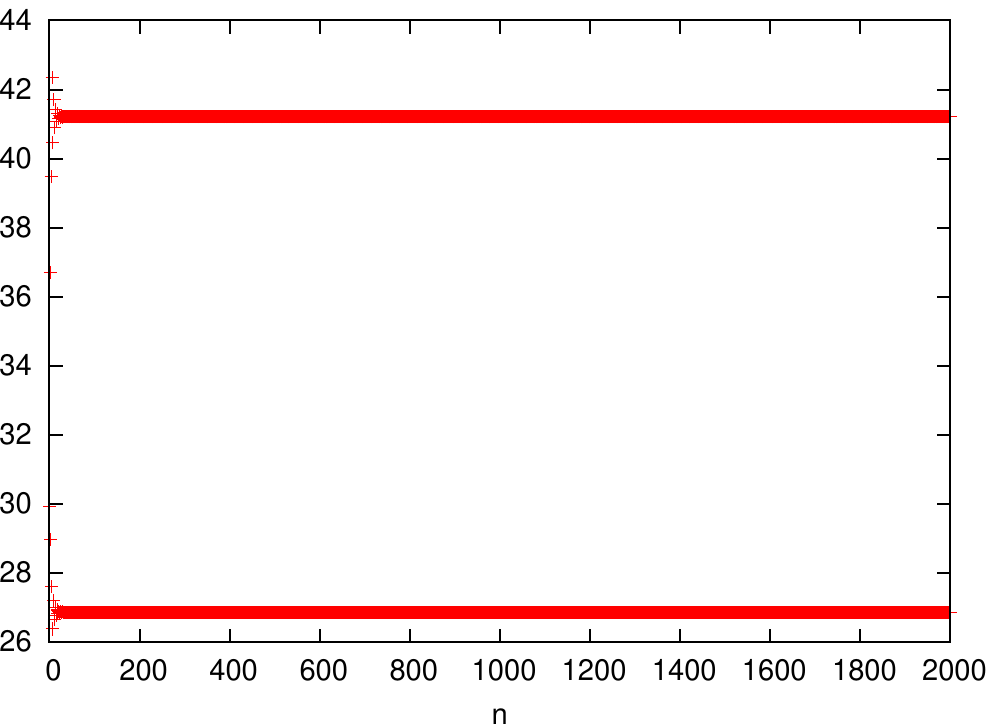}
\includegraphics[height=.23\textheight]{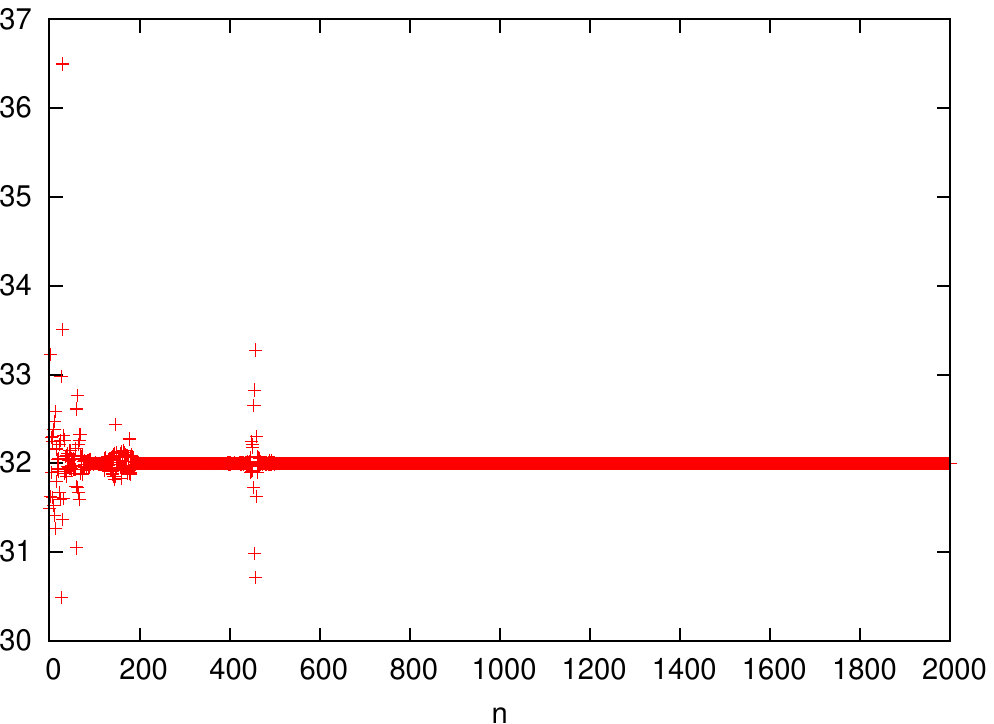}
\caption{Numerical runs for \protect{\eqref{ex:glob}}
in the cases of  (first) no control; (second) $\alpha = 0.417$, $\ell=0$ (third) $\alpha = 0.35$, $\ell=0$
(right) $\alpha = 0.35$, $\ell=0.2$.} 
\label{figure7}
\end{figure}

Fig.~\ref{figure7} (left) illustrates that non-controlled map~\eqref{ex:glob} is chaotic, 
as the theory predicts, there is convergence to the equilibrium without noise for $\alpha=0.417$ (second), and, while for the value of control $\alpha=0.35$ we get a stable two-cycle (third), addition of the Bernoulli noise with $\ell=0.2$ leads to stabilization (right). Note that for $\ell=0.2$, the control value $\alpha=0.35<0.36$ which is theoretically predicted  above.
  
 (b) Consider a continuous function with a unique positive fixed point $K=32$
  \begin{equation}
\label{ex:notglob}
f_2(x):=\left\{\begin{array}{ll}
   50x/28 , & x\in [0, 28),\\
   -5x+190, & x\in [28, 31), \\
   -3x+128, & x\in [31, 32), \\ 
    -2x+96, & x\in [32, 33), \\ 
  -1.7x+86.1, & x\in [33, 45), \\
      9.6, & x\ge 45.
     \end{array}
		\right.
     \end{equation}
     We have the same $\alpha_0=\frac 5{12}\approx 0.416$ as in (a), however $\alpha_0<\Psi(L_2^+, L_2^-) = \Psi(5, 1.7) = 0.4629$. Note that  $\alpha_0$ does not provide global stability: we apply Lemma~\ref{lem:Gull}, and for $x=28$ we  get $G^2(\alpha_0, 28)\approx 26.7455<28$.

\begin{figure}[ht]
\centering
\vspace{-40mm}
\includegraphics[height=.3\textheight]{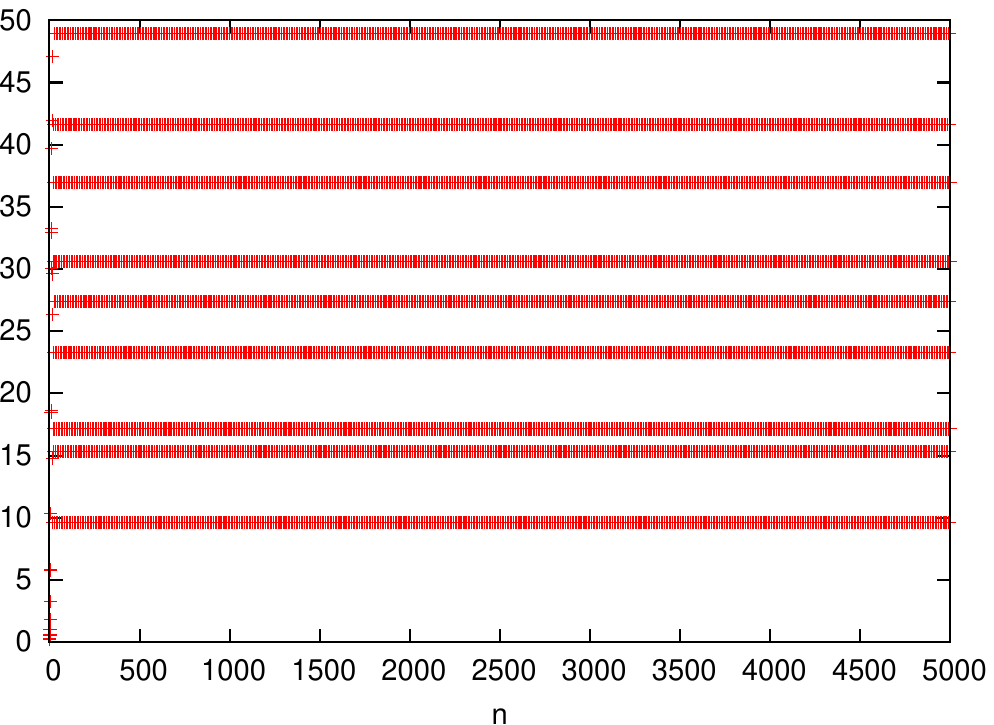}~~~
\includegraphics[height=.3\textheight]{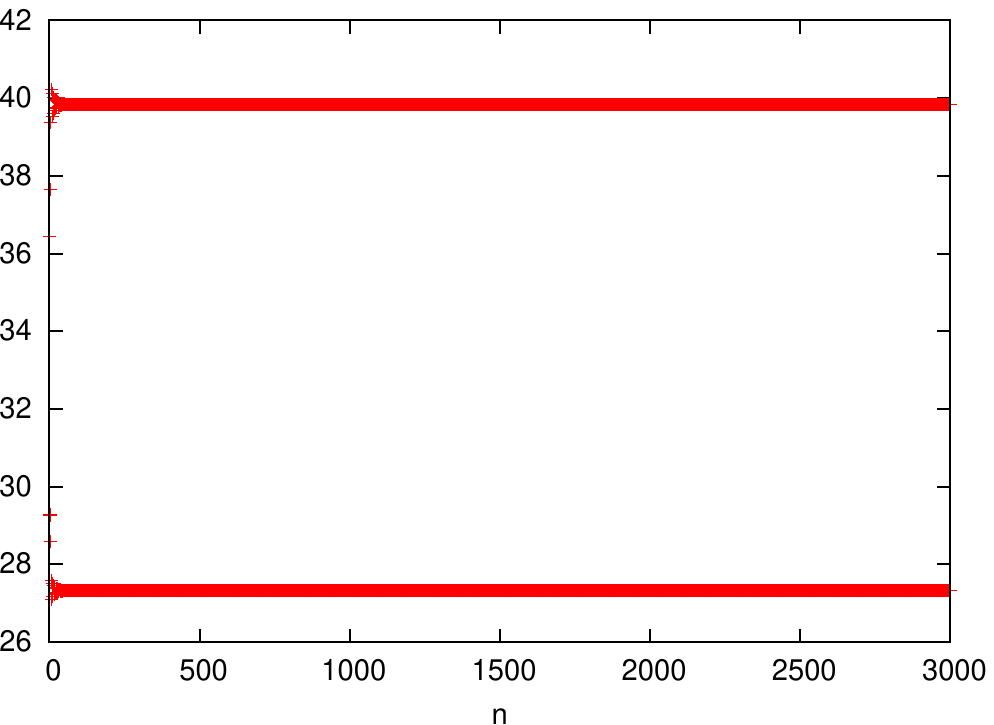}~~~
\includegraphics[height=.3\textheight]{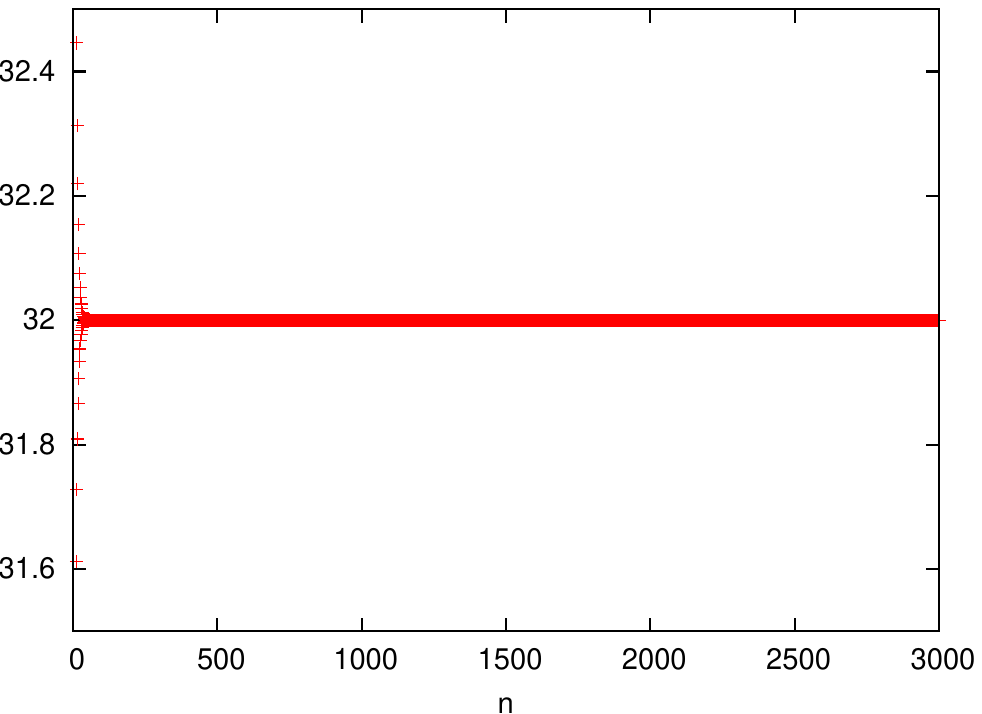}
\caption{Numerical runs for \protect{\eqref{ex:notglob}} 
in the cases when there is no noise and (left) no control; (middle) $\alpha = 0.417$;
(right) $\alpha = 0.463$.}
\label{figure8}
\end{figure}
Fig.~\ref{figure8} (left) illustrates multi-stability and chaotic features of map \eqref{ex:notglob}, the control $\alpha = 0.417$ providing local stability still leads to a two-cycle (middle), while the value of $\alpha = 0.463$ (right) leads to global stability. Fig.~\ref{figure8} illustrates that, generally, local stabilization does not make the equilibrium a global attractor, since
there may also be a stable two-cycle, see also \cite{Singer}.
\end{example}

  \begin{example}
 \label{ex:nonmonot1}
  Consider the function
   \begin{equation}
\label{ex_switching}
   f(x)=\left\{\begin{array}{ll}
	\frac{\pi+2}{\pi-2} x, & x\in (0, 1-2/\pi],\\
   (1-x)\left(1.5+0.5\sin\frac 1{x-1}\right)+1, & x\in (1-2/\pi, 1),\\
   1, & x=1, \\
   (1-x)\left(1.25+0.25\sin\frac 1{x-1}\right)+1, & x\in (1, 1+2/\pi), \\
   1-3/\pi, & x\ge 1+2/\pi.
      \end{array}\right.
  \end{equation}  
  which  is not differentiable at the unique positive equilibrium  $1$, and there is no monotonicity in any neighbourhood of $1$. For its graph see Fig.~\ref{figure_pre_6}.
  We have $L^-=2$, $L^+=1.5$, so, by Proposition~\ref{prop:Globdet},  $\underline\alpha_0=\underline\beta_0=\frac{\tilde L^-\tilde L^+-1}{(\tilde L^-+1)(\tilde L^-+1)}=\frac2{3\times 2.5} \approx 0.266$.
    
In condition \eqref{cond:sides}   we have $a_1= a_2=1$, so in order to apply Theorem \ref{thm:locglgen}, 
we need to have 
    $
    \alpha+\ell<0.5$, $\alpha+\ell> \underline\beta_0=0.266$, $\alpha<\underline\beta_0=0.266$,
  and for the Bernoulli noises,   it should be 
  \begin{equation*}
  \begin{split}
\mathcal V(\alpha, \ell) & := \left[(L^--\alpha(L^-+1))^2-\ell^2(L^-+1)^2  \right]\left[(L^+-\alpha(L^++1))^2-\ell^2(L^++1)^2  \right]\\
   &  = \left[(2-3\alpha)^2-9\ell^2  \right]\left[(1.5-2.5\alpha)^2-6.25\ell^2\right]<1.
  \end{split}
    \end{equation*}
    We can check by  straightforward calculations that the  values 
		$( \alpha, \ell)=(0.23, 0.2)$ and $(0.22, 0.19)$ satisfy the above conditions.
    Thus, introduction of noise decreases the average of the control.

\begin{figure}[ht]
\centering
\vspace{-60mm}
\includegraphics[height=.45\textheight]{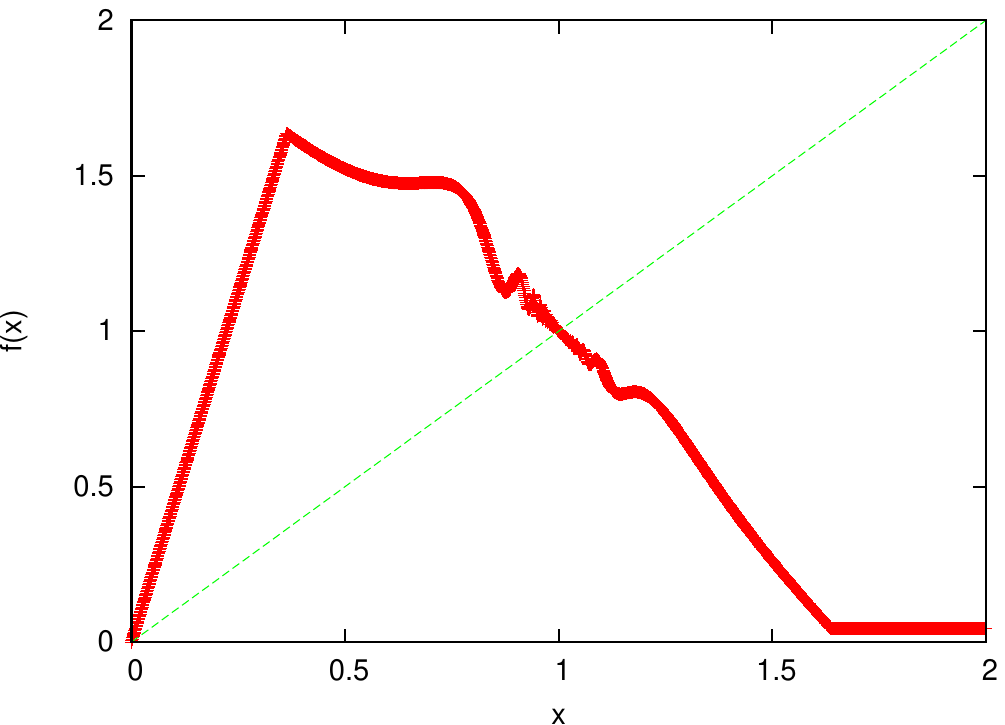}
\caption{The graph of the function 
\protect{\eqref{ex_switching}} without control.}
\label{figure_pre_6}
\end{figure}


\begin{figure}[ht]
\centering
\vspace{-40mm}
\includegraphics[height=.31\textheight]{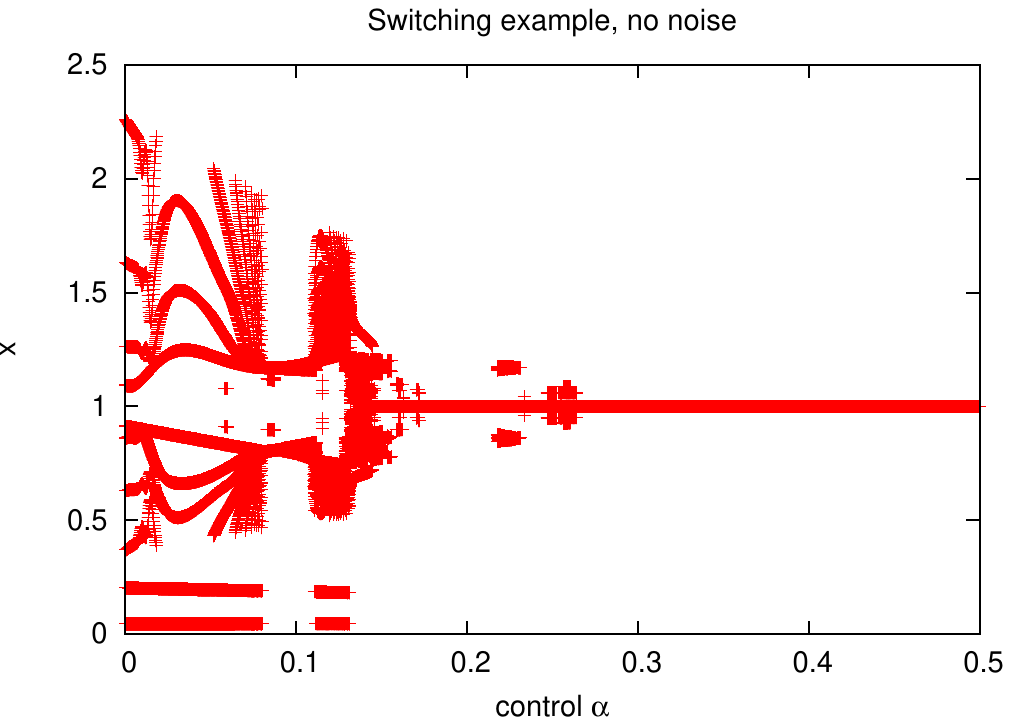}~
\includegraphics[height=.31\textheight]{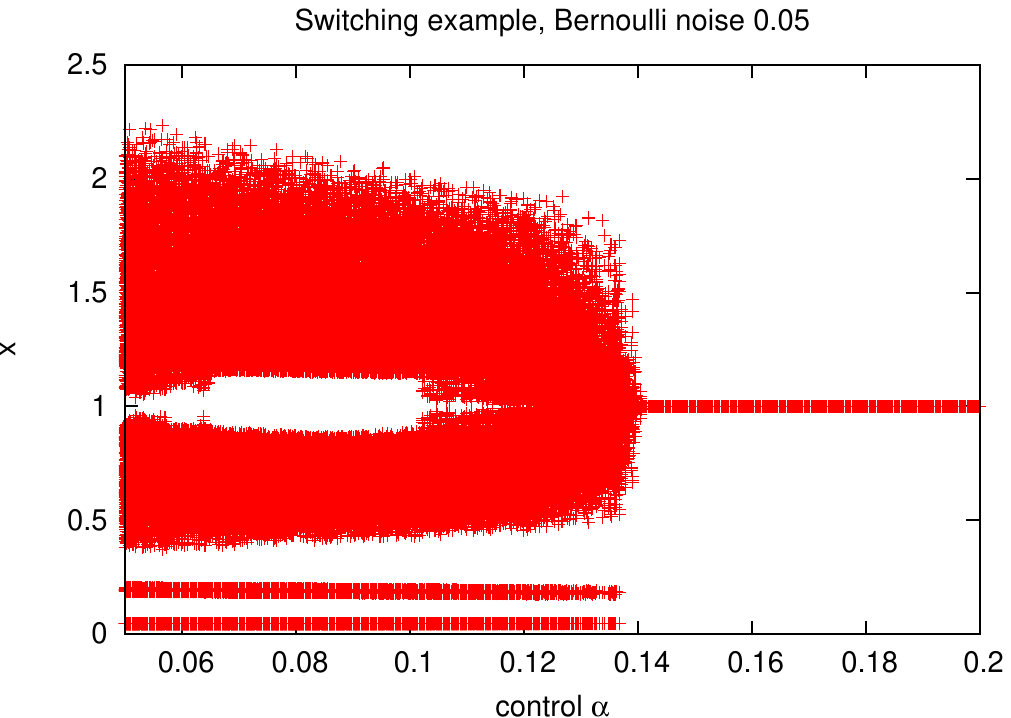}~
\includegraphics[height=.31\textheight]{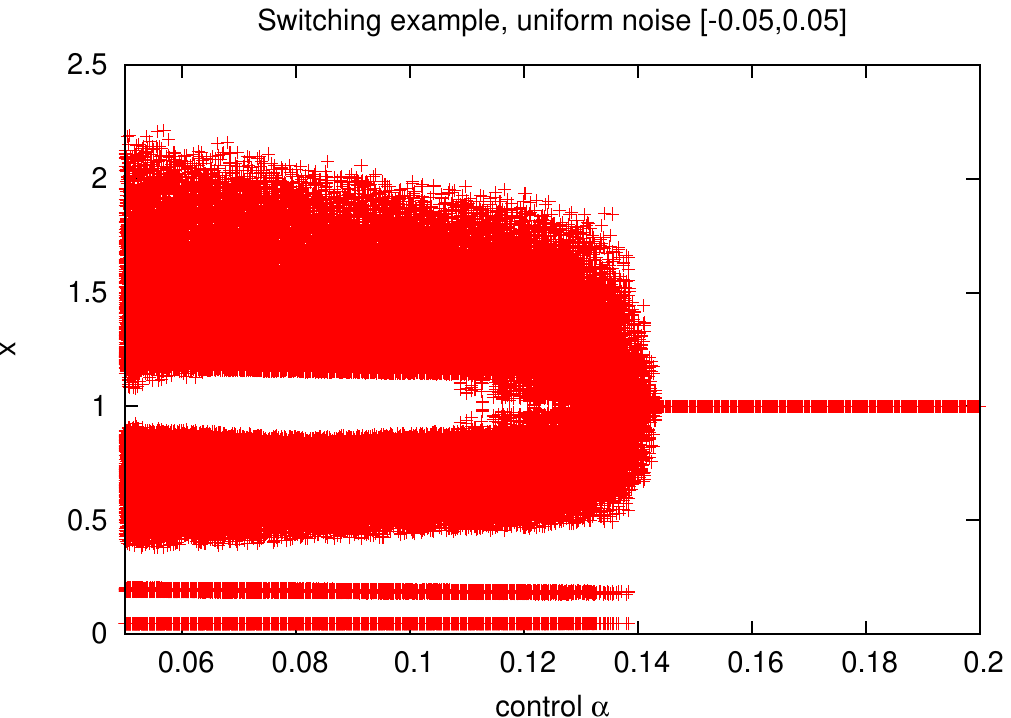}
\caption{Bifurcation diagrams for map \protect{\eqref{ex_switching}} and (left) no noise, (middle) Bernoulli noise taking values $\pm 0.05$, 
(right) continuous noise uniformly distributed on [-0.05,0.05].}
\label{figure6}
\end{figure}
While the bifurcation diagram in Fig.~\ref{figure6} (left) illustrates chaotic behaviour for $\alpha>0.2$, for $\ell=0.05$ stabilization is observed for smaller $\alpha > 0.142$ in the case of the Bernoulli noise, see  Fig.~\ref{figure6} (middle) and for $\alpha > 0.142$ 
in the case of the uniform continuous noise in  Fig.~\ref{figure6} (right). 
These control values are lower than theoretically predicted. 
Possible reasons for this phenomenon are discussed in  Section \ref{sec:disc}.
\end{example}


\section{Conclusions and Discussions}
\label{sec:disc}

Discrete maps are a handy way to describe abundance of semelparous populations. Though simple one-dimensional maps, such as Ricker, for larger values of the parameter can exhibit chaotic behavior, which is not frequently observed in nature. It is sometimes referred to a stabilizing influence of random perturbations, in particular, associated with an environmental noise. In the current paper, we considered control with a prescribed average and  randomness. Can this average be reduced by incorporating noise, and what are the conditions on this average, and admissible noise amplitudes? We give, generally, a positive answer to this question, establishing relevant estimates. This is coherent with experimental observations that introducing noise can either stabilize population size or at least reduce its variation.

The main results of the present paper can be summarized as follows:
\begin{enumerate}
\item
A strict bound for the control parameter leading to global stability is outlined 
in Theorem~\ref{thm:globalvarcontr}, and a smaller average control is allowed in the stochastic settings, provided that $\alpha+ \ell$ is above the critical deterministic bound, see Theorem~\ref{thm:locgldet}.
\item
In the case of a unimodal map, we extend the results of \cite{FL2010,Singer} to the stochastic case 
in Theorem~\ref{thm:Singer}, which is quite a challenging task, taking into account local in general character of convergence in the stochastic case.
To the best of our knowledge, this is the first result of this type. In addition to the general statement, for symmetric distributions, we justify that the average control level ensuring stability is {\em lower in the presence of noise} (Theorem~~\ref{thm:contdiscrdistr}), qualitatively confirming the stabilizing effect of noise.
\item
In the deterministic case, some improvement for global stabilization can be achieved when a series of local Lipschitz constants is taken into account (Proposition~\ref{prop:envG}) in the sense that the control intensity $\alpha$ can be smaller.
\end{enumerate}

The results are verified with numerical simulations illustrating sharpness of the constants when local and global stability are equivalent, and the fact that in the case of sufficient conditions, the required control in examples can be lower than theoretically predicted.
For instance, in Example~\ref{ex:nonmonot1}  the bifurcation diagram in Fig.~\ref{figure6} corresponding to the noise perturbed control demonstrates that stabilization is achieved for  lower than theoretically predicted  average values $\alpha$  of the control.
Note that when there is no noise, the chaotic behavior stops exactly as predicted theoretically, at  $\alpha_0 \approx 0.266$. However,  as soon as noise is present, the necessary average control value is dropped significantly to $\approx 0.142$. We suggest that there could be several  reasons for that. One of them is concerned with the fact that when  the computer does simulations, it assumed that the function  $f$ in \eqref{ex_switching} is equal to 1 in some small, but still significant 
for stochastic stability neighbourhood of 1, so it simulates a slightly different function.  Therefore, in this case,   this is local stability with $\alpha=0$ for the corresponding equation, and for calculation of parameters for the global stability we can apply Theorem \ref{thm:locgldet}, which gives us $\alpha>0.133$ and $\ell>0.266-\alpha$. Taking $\alpha=0.142$  gives $\ell>0.124$. However, the simulation demonstrates that global stability is achieved  for smaller value $\ell=0.05$. We  conjecture that the reason for that is an  oscillating nature of the function \eqref {ex_switching} and the noise lingering at the intervals where the function takes smaller values. 
\medskip
  
The present paper leads to several open questions and lines of research:
\begin{itemize}
\item
For specific types of symmetric distributions, construct sharp stabilization criteria from Theorem~\ref{thm:Singer} and also deduce easily verifiable sufficient conditions. How does the situation change for non-symmetric distributions with the zero mean?
In addition, can the method and the results be generalized to the case of unbounded, for example, normal distributions?
\item
The deterministic results are significantly based on some monotonicity properties of $f$ in some neighbourhood of the unique positive equilibrium. Can the results be adapted to a strongly oscillatory case?   
Example~\ref{ex:nonmonot1} sheds some light on the possibility to extend the results of the present paper to the case when $f$ is oscillatory near $K$.
\item
Everywhere in the current paper we considered the multiplicative noise in the control term. However, it is also interesting as to study additive `environmental' noise 
$$  
x_{n+1}=(1-\alpha)f(x_n)+\alpha x_{n} b+ \ell \xi_{n+1}, \quad x_0>0, \quad \alpha \in [0,1),
$$
for which only a blurred equilibrium can be stabilized in some sense.
\end{itemize}

\section*{Acknowledgment}

The first author acknowledges the support of NSERC, the grant RGPIN-2020-03934.
The authors are grateful to the anonymous referees for their valuable comments.


\section{Apendix}
\subsection{Proof of Lemma~\ref{lem:mathcalL}}

The statements of Parts (i)-(v) are straightforward. By (ii), $\mathcal L^+(1)\mathcal L^-(1)=1$.  Part (v) implies 
$$
\mathcal L^+(\underline \beta_0)\mathcal L^-(\underline \beta_0) = \left[L^+ - (L^+ +1) \left( 1- \frac{1}{1+ L^+} 
- \frac{1}{1+ L^-}  \right) \right] \left[L^- - (L^- +1) \left( 1- \frac{1}{1+ L^+} 
- \frac{1}{1+ L^-}  \right) \right] =1,$$
so (vi) holds. The quadratic polynomial $\mathcal L^+(\beta)\mathcal L^-(\beta)$ has a positive coefficient of $\beta^2$ and two roots: one and $\beta_0<1$. Thus,
$\mathcal L^+(\beta)\mathcal L^-(\beta)<1 $ holds when $\beta\in (\underline \beta_0, 1)$, which concludes the proof of Part (vii).

\subsection{Proof of Lemma~\ref{lem:aux21}}
If $x_0\in [f^2_m,  f_m]$ the result follows from Lemma~\ref{lem:Gv}~(iv) with $S_0=1$.  
Assume that $x_0\in (0, f^2_m)$ and 
let $\Delta_1:=\inf\{f(x)-x ~:~ x\in [x_0, f^2_m) \}>0$ as $f^2_m<K$, 
$G(\alpha_n, x)-x=(1-\alpha_n)(f(x)-x)$, for each $\alpha_n\in (\beta_*, \beta^*)$ and $1-\alpha_n>1-\beta^*$. Then
\[
G(\alpha_n, x)-x>(1-\beta^*)\Delta_1, \, \, x\in [x_0, f^2_m], \quad N^-:=\left[\frac{f^2_m-x_0}{\Delta_1(1-\beta^*)}\right]  +1.
\]
Reasoning inductively and assuming that $x_i=G(\beta_{i}, x_{i-1})\in (x_{i-1}, f^2_m]$, $i=1, \dots, k$, $k\le N^-$, we get 
\begin{align*}
&x_{k+1}=G(\alpha_{k+1}, x_{k})-x_{k}+x_{k}\ge  G(\beta^*, x_{k})-x_{k}+x_{k}\ge\Delta_{1}+x_{k}\ge \dots \ge k\Delta_{1}+x_{0}. 
\end{align*}
So after at most $N^-$ steps the solution reaches $ [f^2_m,  f_m]$.  Recall that by definition of $f_m$ the solution cannot jump over $f_m$.
In the case $x_0>f_m$  we find $N^+$  in  a similar way and let $S_0:=\max\{N^-, N^+\}$.  

\subsection{Proof of Lemma~\ref{lem:Gullbeta}}
Note that it is enough to consider only $x\in [f_m^2, f_m]$. 
Fix some $\alpha>\frac{L^-}{L^-+1}$, then, by Lemma~\ref{lem:mathcalL}~(i),(iv) we have  $\mathcal L^+(\alpha)<\mathcal L^-(\alpha)<0$, and  
\begin{align*}
G(\alpha, x)-K=(1-\alpha)(f(x)-K)+\alpha(x-K)<[(1-\alpha)L^--\alpha](K-x)=\mathcal L^-(\alpha)(K-x)<0,
\, x\in (f^2_m, K),\\
K-G(\alpha, x)<[(1-\alpha)L^+-\alpha](x-K)=\mathcal L^+(\alpha)(K-x)<0,\quad x\in (K, f_m).
\end{align*}
This  implies $ x<G(\alpha, x)<G^2(\alpha, x)$ for $x\in (f^2_m, K)$ and $ x>G(\alpha, x)>G^2(\alpha, x)$ for $x\in (K, f_m)$, so $\alpha\in \mathcal S$. 

Further, we notice that, since $\mathcal S\neq \emptyset$, 
for each $\alpha>\underline \beta_*= \inf \mathcal S$ there is $\beta_1\in \mathcal S$ s.t. $\alpha>\beta_1$, which, by Lemma~\ref{lem:Gv}(ii), implies that $G(\alpha, x)>G(\beta_1, x)$ for  $x>K$ and  $G(\alpha, x)<G(\beta_1, x)$ for  $x<K$.  If $x\in (f^2_m, K)$ and  $G(\alpha, x)>K$, there is $\hat x\in (x, K)$, s.t. $G(\alpha, x)=G(\beta_1, \hat x)$. Since $\beta_1\in \mathcal S$ we have  $G^2(\alpha, x)=G(\alpha, G(\beta_1, \hat x))>G^2(\beta_1, \hat x)>\hat x>x$, which, by Lemma~\ref{lem:Gull}, yields  the result. Other cases are either similar or have been considered above. 
Note that the proof of the second part of Lemma~\ref{lem:Gullbeta} follows the scheme of \cite[Theorem 3]{Gull}, even though we do not impose any monotonicity restriction on $f$. 

\subsection{ Proof of Lemma~\ref {lem:auxback}}
Once $x_n-K$ does not change the sign, or if it changes the sign once, it is true. Assume $x_n<K$, the case $x_n>K$ is similar.
For the proof it is enough to show that  when a solution moves to the right of $K$ and then to the left of $K$, at the moment of the first return  to $(0, K)$ it will be on the right of $x_n$.
Let $s_1=\min\{s>n: x_{s}>K\}$, $s_2=\min\{s>s_1: x_{s}<K\}$, and both sets be not-empty. Let $s_1>n+1$, $s_2>s_1+1$. Then, 
   $
   x_n<x_{n+1}\le x_{s_1-1}<K<x_{s_1}, \quad x_{s_1}>x_{s_1+1}\ge x_{s_2-1}>K>x_{s_2}.
   $   
   Since  $x_{s_1}>x_{s_2-1}>K>x_{s_1-1}$  and $\alpha_{s_2}>\beta_*$,  we have $x_{s_1}=G(\alpha_{s_1}, x_{s_1-1})\le G(\beta_*, x_{s_1-1})$. Due to the continuity of $G(\beta_*, \cdot)$,  
    there is $\bar x\in (x_{s_1-1}, K)$ s.t. $G(\beta_*, \bar x)=x_{s_2-1}$. Since $\beta_*> \underline\beta_*$ we have $G^2(\beta_*, \bar x) >\bar x$ and since $x_{s_2-1}>K$ we have $G(\alpha_{s_2}, x_{s_2-1})\ge G(\beta_*, x_{s_2-1})$.  Then
 $ x_{s_2}=G(\alpha_{s_2}, x_{s_2-1})\ge G(\beta_*, x_{s_2-1})=G(\beta_*, G(\beta_*, \bar x))=G^2(\beta_*, \bar x) >\bar x>x_{s_1-1}>x_n$,
 so $x_{s_2}>x_n$, i.e.  the first return of the solution to $(0, K)$ is to the right of $x_n$.
 
 If $s_1=n+1$,  $s_2>s_1+1$, we have $x_n<K<x_{n+1}$, and we can start the reasoning  as in the first case for  $x_{n+1}$ instead of $x_n$.   If $s_1=n+1$, $s_2=s_1+1$,  we have $x_n<K<x_{n+1}, \quad x_{n+2}<K<x_{n+1}$. So  $x_{n+1}=G(\alpha_{n+1}, x_{n})\le G(\beta_*, x_{n})$ and for some $\bar x\in (x_{n}, K)$ we have $G(\beta_*, \bar x)=x_{n+1}$.  Thus
$
 x_{n+2}= G(\alpha_{n+2}, x_{n+1})\ge G(\beta^*, x_{n+1})=G(\beta^*,G(\beta_*, \bar x))> \bar x >x_n,
  $
 and again  the first return of solution to $(0, K)$ is to the right of $x_n$. 

 
 \subsection{ Proof of Theorem \ref{thm:globalvarcontr}}
 
 
   (i) Consider first $x_0\in [f^2_m, f_m]$. Assume the contrary that a solution does not converges to $K$.
		By Assumption~\ref{as:Lglob1}, there is an infinite number of points $x_n$ on both sides of $K$,  $(x_n)_{n\in {\mathbb N}_0}=(x_{n_i})_{i\in {\mathbb N}_0}\cup (x_{n_j})_{j\in {\mathbb N}_0}$, $x_{n_i}\in (0,K)$, $x_{n_j} \in (K, \infty)$.    
Lemma~\ref{lem:auxback} implies that a solution   
to \eqref{eq:PBCGstoch} is a union of two monotone subsequences: the first one is in $(0,K)$ and increasing, the second  one is  in $(K, \infty)$ and decreasing. 
        If  $(x_n)_{n\in \mathbf N}$ does not converge, there are non-negative $a$ and $b$ and $\bar N\in \mathbb N$,    s.t.  
 $x_n\notin [K-a, K+b]$ for all $n\ge \bar N$. If we assume that one of $a,b$ is equal to zero, we get that one of the sequences, say $x_{n_i}$, converges to $K$, then from continuity of $f$, the other sequence of  $f^m(x_{n_i})$
also converges to $K$.

Define
    \begin{equation}
    \label{def:Delta*}
    \Delta_*:=\min\left\{\inf_{u\in  [f^2_m, K-a]\cup [K+b, f_m]} \left| G(\beta_*, u)-u \right|, \quad \inf_{u\in [f^2_m, K-b/{L^-}]\cup [K+a/{L^-}, f_m]} \left| G^2(\beta_*, u)-u  \right| \right\},
    \end{equation} 
  and note that $\Delta_*>0$  by Assumptions \ref{as:Lglob1}, Lemma~\ref{lem:Gull}, definition \eqref{def:beta*Gull}, and by the choice of $\beta_*$.
  
  Let $x_i\in  [f^2_m, K-a]$ for some $i> \bar  N$.  
  Applying Lemma~\ref{lem:Gv}~(iii), we get
    \begin{align*}
x_{i+1}-x_i = G(\alpha_{i+1}, x_i)-x_i & =\left(1-\frac{\alpha_{i+1}-\beta_*}{1-\beta_*}\right)G(\beta_0, x_i)+\frac{\alpha_{i+1}-\beta_*}{1-\beta_*} x_i-x_i\\ & =
\frac{1-\alpha_{i+1}}{1-\beta_*}[G(\beta_*, x_i)-x_i]>\frac{1-\beta^*}{1-\beta_*}\Delta_{*}=\Delta_{*}.
  \end{align*}  
 If in addition, $G(\alpha_{i+1},x_{i})>K+b$, there is $\bar x_i\in (x_i, K)$ s.t. $G(\alpha_{i+1},x_{i})=G(\beta_*, \bar x_i)$.   By Assumption \ref{as:Lglob2}  we have
    $
    G(\beta_*, \bar x_i)\le K+L^-(K-\bar x_i)$, so $K+b< G(\alpha_{i+1},x_{i})=G(\beta_*, \bar x_i)\le K+L^-(K-\bar x_i),
    $
  which implies 
$$
  b<L^-(K-\bar x_i), \quad K-\bar x_i>\frac b{L^-}, \quad \bar x_i\in \left(x_i, K-\frac b{L^-} \right)\subset \left(f^2_m, K-\frac b{L^-} \right). 
$$
    Recalling  that $G(\alpha_i, x)>G(\beta_*, x)$ for $x>K,$ we get 
 $$
 G(\alpha_{i+1}, G(\alpha_i, x_i))>G(\beta_*, G(\alpha_i, x_i))=G(\beta_*, G(\beta_*,\bar x_i))=G^2(\beta_*,\bar x_i)>\bar x_i>x_i,
 $$
which, by \eqref{def:Delta*}, implies 
$
 G(\alpha_{i+1}, G(\alpha_i, x_i))-x_i>G^2(\beta_*,\bar x_i)-\bar x_i>\Delta_*.$ So, if the solution changes the side of $K$ at two successive steps, returning to $[f^2_m, K-a]$, we have $x_{i+2}>x_i+\Delta^*$. Following the same argument as in the proof of  Lemma~\ref{lem:auxback}, we actually can  get $x_{i+m}>x_i+\Delta^*$, where $i+m$  is the first moment after returning to $[f^2_m, K-a]$.
 
  Therefore  the solution $x_i$, which starts in $[f^2, K-a]$, moves right with the step of the length bounded below by $\Delta_*$. If it remains on the left of $K$, it eventually  moves to the right of $K-a$ in a finite number of steps, which contradicts to our assumption. If, at some moment it jumps over $K+b$ and remains to the right of $K$, it eventually  moves below $K+b$ in a finite number of steps, which contradicts to our assumption again. If it returns to $[f^2, K-a]$, its new position there will be at least $\Delta_*$ to the right than before the jump. Summarizing all the above, we conclude that the solution cannot move in $[f^2, K-a]$ for more than $N_1:=\frac{K-a-f^2_m}{\Delta_*}$  steps, cannot stay in $[K-b, f_m]$ for more than $N_2:=\frac{f_m-K-b}{\Delta_*}$  steps, and cannot remain in  $[f^2_m, K-a]\cup [K+b, f_m]$ for more than $N_1+N_2+\frac{f_m+a-f^2_m-b}{\Delta_*}$  steps.   
 
 Similarly, when $x_i\in  [K+b, f_m]$  we get $x_i-x_{i+1}>\Delta_* $.    If in addition, $G(\alpha_{i+1},x_{i})<K-a$, there is $\bar x_i\in (K, x_i)$ s.t. $G(\alpha_{i+1},x_{i})=G(\beta_*, \bar x_i)$.   By Assumption \ref{as:Lglob2},  we get
$K-a> G(\alpha_{i+1},x_{i})=G(\beta_*, \bar x_i)\ge K-L^+(\bar x_i-K)$, so $\bar x_i\in \left( K+\frac a{L^+},  x_i\right)\subset \left(K+\frac a{L^+}, f_m \right),$  and
$
x_i- G(\alpha_{i+1}, G(\alpha_i, x_i))>\bar x_i-G^2(\beta_*,\bar x_i)>\Delta_*.$  
 
The case $x_0\in (0, f^2_m)\cup (f_m, \infty)$  follows from Lemma~\ref{lem:aux21}.
   
 (ii) 	Basically we repeat the proof of Part (i) for $a=b=\delta$. The necessary number of steps in this case is $\max\{\bar N_1, \bar N_2\}$, where $\bar N_1$ is the first moment when the solution gets into $[K-\delta, K)$, and $\bar N_2$ is the first moment when the solution gets into $(K,  K+\delta)$, if     both numbers $\bar N_1$, $\bar N_2$ are finite. If one of them is infinite, say $\bar N_2=\infty$, we put $S_1=\bar N_1$. In other words,
  $S_1:=\max\{\bar N_1, \bar N_2\}$ if $\bar N_2, \bar N_1<\infty$,  $S_1=\bar N_1$  if $\bar N_{2}=\infty$, and $S_1=\bar N_2$ if $\bar N_{1}=\infty$.

  \subsection{Proof of Theorem~\ref {thm:locgldet}}
  For any  $(\alpha, \ell)$ satisfying \eqref{cond:ael} we have  $0<\alpha_0<\alpha-\ell<\alpha_n<\alpha+\ell<1$, for all $n\in \mathbb N$,  and $\frac{\underline\beta_*-\alpha}{\ell}<1.$
    We consider only $\alpha<\underline\beta_*$, since otherwise  the solution  to \eqref{eq:PBCintr}   is globally asymptotically  stable for $\ell=0$. Then, by Assumption~\ref{as:noise},
 $ \mathbb P(\hat \Omega)>0$, where $\hat \Omega=\left\{\omega\in \Omega:\xi(\omega)\in \left( \frac{\underline\beta_*-\alpha}{\ell},  1\right]\right\}$,
    and by Lemma~\ref{lem:topor} there is a random moment $\mathcal N$ s.t. 
    \begin{equation}
    \label{def:momentN}
    \xi_{\mathcal N+S_0}, \xi_{\mathcal N+1+S_0},  \dots,  \xi_{\mathcal N+S_1+S_0}  \in  \left( \frac{\underline\beta_*-\alpha}{\ell},  1\right],
     \end{equation}	
     where $S_0$ and $S_1$ are from Lemma~\ref{lem:aux21} and Theorem~\ref{thm:globalvarcontr}~(ii), respectively, $S_0=S_0(\alpha-\ell,  \alpha+\ell, x_0)$, $S_1=S_1(x_0, \alpha-\ell, \delta_0)$.
     Fix some  $k\in \mathbb N$, set  
     \begin{equation}
  \label{def:Omegak}    
   \Omega_k=\left\{\omega\in \Omega:   \mathcal N=k \right\}=\left\{\omega\in \Omega:    \xi_{k+i+S_0}  \in  \left( \frac{\underline\beta_*-\alpha}{\ell},  1\right],\, i=0, \dots, S_1\right\},
     \end{equation}
      and note that $\Omega_k$ is defined by $\xi_{k+S_0}, \xi_{k+1+S_0},  \dots,  \xi_{k+S_1+S_0}$.
      
      By Lemma~\ref{lem:aux21}, $x_{S_0+k}\in[f^2_m, f_m]$  on all $\Omega$. 
   Let $y$ be a solution to 
   \[
   y_{n+1}=G(\bar \alpha_n, y_n), \quad y_0=x_{S_0+k}, \quad \bar \alpha_n=\alpha+\ell \xi_{S_0+k+n},
   \]
   considered path-wise on $\Omega_k$. 
	Since $\alpha_n=\alpha+\ell \xi_n>\underline\beta_*$ for $n=S_0+k+1, S_0+k+2, \dots,S_0+k+S_1$, \, 
Theorem~\ref{thm:globalvarcontr}~(ii) implies  $x_{S_0+k+S_1}=y_{S_1}\in  (K-\theta, K+\theta)$ on   $\Omega_k$.   Since $\alpha>\alpha_0$, as soon as  $x$ gets into $(K-\theta, K+\theta)$,  it tends to $K$. 
So, for each  $\omega\in\Omega_k$ we get $\displaystyle \lim\limits_{n\to \infty} x_n=K$. Since $\Omega=\cup _{k=1}^\infty\Omega_k$, this completes the proof.
\subsection{Proof of Lemma~\ref{lem:locdet}}
As we have assumed everywhere,  $\tilde L^-\ge \tilde L^+$ and by Lemma~\ref{lem:mathcalL}, $\tilde {\mathcal L}^+(\alpha)<1$, $\tilde {\mathcal L}^-(\alpha)\tilde{\mathcal L}^+(\alpha)<1$. In order to keep solution inside of $(K-\theta, K+\theta)$, where inequalities \eqref{cond:loc2} can be applied, we need to decrease the left part of the interval. So, if $x\in (K-\theta/ \tilde{\mathcal L}^-(\alpha), K)$ and $G(\alpha, K)>K$, we have
$
G(\alpha, K)-K\le \tilde {\mathcal L}^-(\alpha)(K-x)<\tilde {\mathcal L}^-(\alpha) \, \, \theta/\tilde{\mathcal L}^-(\alpha)
=\theta,
$
i.e. $G(\alpha, K)\le K+\theta.$

\subsection{Proof of Theorem \ref{thm:locglgen}} 
  Suppose that the statement of theorem  does not hold, i.e. there exists a pair $(\alpha, \ell)$ satisfying  \eqref{cond:lambda0}, \eqref{cond:sides}, and \eqref{cond:glob}  such that 
     \begin{equation}
 \label{cond:notglob} 
 \begin{split}
  \mbox{for some  $x_0>0$, $\kappa\in (0, 1)$, $\Omega_\kappa\subset \Omega$ with $\mathbb P (\Omega_\kappa)=\kappa$,  the solution $x_n \not \to K$ on $\Omega_\kappa$.}
  \end{split}
  \end{equation}
  Without loss of generality we can assume that $\kappa\in (0, 2/3)$.
  In the proof below we consider a solution to \eqref{eq:PBCintr} with  $(\alpha, \ell)$ 
	and $x_0$ satisfying  \eqref{cond:notglob}.
  Note that, once $x_n=K$, all $x_{n+j}=K$, $j \in {\mathbb N}$, so we only have to consider the case $x_n \neq K$,
	$n \in {\mathbb N}$.
    
  (i) We start with the proof that a solution to \eqref{eq:PBCintr}  changes sides of $K$  at each step.
  We have 
     \begin{align*}
     &G(\beta, x)-K=(1-\beta)(f(x)-K)+\beta(x-K)\ge [(1-\beta)a_1-\beta](K-x), \quad \mbox{if} \quad x\in [K-\theta, K],  \nonumber\\   
     &K-G(\beta, x)=(1-\beta)(K-f(x))+\beta(x-K)\ge [(1-\beta)a_2-\beta](x-K), \quad \mbox{if} \quad x\in [K, K-\theta].
     \end{align*}
     For $\beta=\alpha+\ell \xi$ we get $\alpha-\ell\le \beta\le\alpha+\ell$, so
     \[
     (1-\beta)a_1-\beta=a_1-\beta(a_1+1)\ge a_1-(\alpha+\ell)(a_1+1), \quad   (1-\beta)a_2-\beta\ge a_2-(\alpha+\ell)(a_2+1),
        \]
which,  along with \eqref{cond:sides}, implies $(G(\alpha+\ell \xi_n, \, x)-K)(K-x)>0$, $x \in (K-\theta, K)\cup(K, K+\theta)$, $n\in \mathbb N$. Therefore, as soon as a solution remains in $(K-\theta, K+\theta)$, it changes position relative to $K$ at each step.  
Since 
\begin{equation*}
\begin{split}
&G(\alpha+\ell \xi_n, \, x)-K\le {\mathcal L}^-(\alpha+\ell \xi_{n})(K-x), \quad x\in (K-\theta, K),\,\\
&K-G(\alpha+\ell \xi_n, \, x)\le {\mathcal L}^+(\alpha+\ell \xi_{n})(x-K), \quad x\in (K, K+\theta),
\end{split}
\end{equation*}
we conclude that ${\mathcal L}^-(\alpha+\ell \xi_{i})>0$ and ${\mathcal L}^+(\alpha+\ell \xi_{i})>0$,
so  we can omit the absolute value sign in \eqref{cond:lambda0}.

(ii) Now, let us prove local stability. 
Consider the  sequence  $(u_i)_{i\in \mathbb N}$ of i.i.d. variables
  \begin{equation}
  \label{def:ui}  
   u_i:=\ln \left[{\mathcal L}^-(\alpha+\ell \xi_{i}) {\mathcal L}^+(\alpha+\ell \xi_{i+1})\right], \quad  i\in \mathbb N.
 \end{equation}
By monotonicity of $\mathcal L$, see Lemma~\ref{lem:mathcalL}~(i), we have ${\mathcal L}^-(\alpha+\ell \xi_{t+i})\le {\mathcal L}^-(\alpha-\ell)$, ${\mathcal L}^+(\alpha+\ell \xi_{t+i})\le {\mathcal L}^+(\alpha-\ell),$
for any $i\in \mathbb N$.
 By \eqref{cond:lambda0} we have $\mathbb E u_i=-\lambda_0$. Based on Corollary \ref{cor:Kolm}, for each $\nu>0$ we can find a nonrandom number $\bar N=\bar N(\nu)$ such that, $\mathbb P \left\{\Omega_{\nu, 0}\right\}> 1-\nu$, where 
$\Omega_{\nu, 0}:=\left\{\omega \in \Omega: \sum_{i=1}^{n}  u_i\le  -\lambda_0 \frac{n}{2} \quad \text{for all} \quad  n\ge \bar N  \right\}.$
For $t\in \mathbb N$, set
\begin{equation}
\label{def:Omeganut}
\Omega_{\nu, t}:=\left\{\omega \in \Omega: \sum_{i=t}^{t+n-1}  u_i\le  -\lambda_0 \frac{n}{2} \quad \text{for all} \quad  n\ge \bar N  \right\}.
\end{equation}
In general, $\Omega_{\nu, 0}\neq \Omega_{\nu, t}$ for $t\neq 0$, but since $u_i$ are identically distributed, we have  $\mathbb P(\Omega_{\nu, 0})= \mathbb P(\Omega_{\nu, t})>1-\nu$ for  each $t\in \mathbb N$. Also, 
$\displaystyle e^{\sum_{i=t}^{t+n-1}  u_i} \le  e^{-\lambda_0 n/2}$, for all $n\ge \bar N$, on  $\Omega_{\nu, t}$.
For $\kappa$ from \eqref{cond:notglob}, we set 
 \begin{equation}
 \label{def:delta0} 
 \nu:=\frac{\kappa}4, \quad 
  \delta_0:=\theta \mathcal B^{-\bar N}, \, \mbox{where}\,\,
 \mathcal B:=\max\{{\mathcal L}^-(\alpha-\ell), \, {\mathcal L}^+(\alpha-\ell)\}.
 \end{equation}
Let us demonstrate that, as soon as $x_t\in (K-\delta_0, K+\delta_0)$, we get $\lim\limits_{n\to \infty} x_n=K$ on $\Omega_{\nu, t}$. 
 Assume for simplicity that  $\bar N$ and $s$ are even, $\bar N=2\bar M$, $s=2d$.   
 We have  $\Omega=\Omega^+(t)\cup \Omega^-(t)$, where 
  $\Omega^+(t)=\{\omega\in \Omega:x_t(\omega)\in (K, K+\delta_0) \}$, $ \Omega^-(t)=\{\omega\in \Omega:x_t(\omega)\in (K-\delta_0, K) \}.$
On  $\Omega^-(t)\cap \Omega_{\nu, t}$  we have $x_{t+1}>K$, $x_{t+2}<K$, {\it etc}, if the solution remains in $(K-\delta_0, K+\delta_0)$, so, 
inductively,
 \begin{equation*}
 \begin{split}
 \label{est:loc1}
 x_{t+1}-K & \le {\mathcal L}^-(\alpha+\ell \xi_{t+1})[K-x_t]\le {\mathcal L}^-(\alpha-\ell)[K-x_t]<\mathcal B\delta_0<\theta, \\
 K-x_{t+2} & \le {\mathcal L}^+(\alpha+\ell \xi_{t+2})|x_{t+1}-K|< {\mathcal L}^+(\alpha+\ell \xi_{t+2}) {\mathcal L}
^-(\alpha+\ell \xi_{t+1})|x_t-K|\\&=e^{u_t}|x_{t}-K| <\mathcal B^2\delta_0<\theta, \\
&\dots \dots \dots \dots \dots \dots\\
|x_{t+\bar N}-K| & \le e^{\sum_{i=t}^{t+\bar M-1}u_i}|x_t-K|\le \mathcal B^{\bar N}\delta_0=\theta.
 \end{split}
  \end{equation*}

 Using  \eqref{def:ui} and continuing estimations, we arrive at 
    \begin{equation*}
 \label{est:loc11} 
|x_{t+\bar N+d}-K|= |x_{t+2(\bar M+s)}-K|
 \le \exp \left\{ \sum_{i=t}^{t+\bar M+s-1}u_{i} \right\} |x_t-K|<e^{-\lambda_0 (\bar M+s)/2}\delta_0<\theta. 
  \end{equation*}
 Similar inequalities can be obtained on $\Omega^+(t)\cap \Omega_{\nu, t}$. So,   $
      |x_{t+n}-K|\le   e^{-\lambda_0 (n)/2}\delta_0\to 0$, as $n\to \infty$,  on $\Omega_{\nu, t}$. 

(iii) Now proceed to the proof of global attractivity. 
Let $S_0=S_0(\alpha-\ell, \alpha+\ell, x_0)$, $S_1=S_1(x_0, \alpha-\ell, \delta_0)$ be from Lemma~\ref{lem:aux21} and Theorem~\ref{thm:globalvarcontr} (ii), respectively.  Note  that $\delta_0$ was chosen as in \eqref{def:delta0}, so $\delta_0$ and $S_1(x_0, \alpha-\ell, \delta_0)$ depend on $\kappa>0$, which is the lower estimate for the probability  of the set $\Omega_\kappa$, where $x_n \not \to K$.

 Recall that $0<\alpha-\ell<\alpha_n<\alpha+\ell<1$, for all $n\in \mathbb N$,  and $\frac{ \underline\beta_*-\alpha}{\ell}<1.$ 
 Reasoning as in the proof of Theorem \ref{thm:locgldet} and using the same notations  \eqref{def:momentN} and \eqref{def:Omegak} for the random moment $\mathcal N$ and  sets $\Omega_k$, respectively,   we conclude    that $x_{S_0+j+S_1}\in  (K-\delta_0, K+\delta_0)$, on   $\Omega_j$.   Denoting $t=t(j):=S_0+j+S_1,$
   and considering $\nu$ and $\Omega_{\nu, t}$ defined as in  \eqref{def:delta0} and \eqref{def:Omeganut}, respectively,  with $\mathbb P(\Omega_{\nu, t})>1-\nu$, we arrive at 
 $\lim\limits_{n\to \infty} x_n=K$,  on $\Omega_{\nu, t}\cap \Omega_k$.
 Since  $\Omega_{\nu, t}$ is defined by $\{\xi_i, i > t(j)=S_0+j+S_1\}$, while $\Omega_k$ consists of $\{\xi_i, i\le t(j)=S_0+j+S_1\}$, by independence of $\xi_i$ and by definition \eqref{def:delta0}, we have 
   $
  \mathbb P( \Omega_{\nu, t(k)}\cap \Omega_k)= \mathbb P( \Omega_{\nu, t(k)})  \mathbb P(\Omega_k)\ge
 (1-\nu) \mathbb P(\Omega_k)=(1-\kappa/4)\mathbb P(\Omega_k).
     $
    Since $\cup_{j=1}^\infty \Omega_j=\Omega$, we can choose $j_\kappa$ s.t. $
    \cup_{j=1}^{j_\kappa} \mathbb P\left(\Omega_j \right)>\frac{1-\kappa/2}{1-\kappa/4}$. Letting $\tilde \Omega:=\cup_{j=1}^{j_\kappa} \left[ \Omega_{\nu, t(j)}\cap \Omega_j\right]$, we arrive at 
\[
 \mathbb P(\tilde \Omega)= \mathbb P\left(\cup_{k=1}^{k_\kappa} \left[ \Omega_{\nu, t(k)}\cap \Omega_k\right]\right)=\sum_{k=1}^{k_\kappa} \mathbb P\left(\Omega_{\nu, t(k)}\cap \Omega_k\right)\ge (1-\kappa/4)\sum_{k=1}^{k_{\kappa}}  \mathbb P\left(\Omega_k\right)>1-\frac{\kappa}2
 \]
 and  $\lim\limits_{n\to \infty} x_n=K$,  on $\tilde\Omega$. However, by our assumption in \eqref{cond:notglob}, we should have 
 $\tilde \Omega\subset \Omega\setminus \Omega_\kappa$,  so $1-\frac{\kappa}2\le \mathbb P(\tilde \Omega)\le \mathbb P(\Omega\setminus \Omega_\kappa)=1-\kappa$. The contradiction proves that $\lim\limits_{n\to \infty} x_n=K$ a.s.
\subsection{Proof of Theorem~\ref{thm:Singer}}
Since $f$ is  continuously differentiable at $K$ and $L_0>1$, for each $\varepsilon\in (0, 1)$ 
there exists $\theta=\theta(\varepsilon)>0$  such that
     \begin{equation}
     \label{cond:Lipscheps}
     \begin{split}
& \left(L_0-\varepsilon\right)(K-x)<f(x)-K<\bigl(L_0+\varepsilon\bigr)(K-x), \quad  x\in (K-\theta, K), \\
& \bigl(L_0-\varepsilon\bigr)(x-K)< K-f(x)<\bigl(L_0+\varepsilon\bigr)(x-K), \quad  x\in (K, K+\theta).
      \end{split}
      \end{equation} 
Relations \eqref{cond:Lipscheps} imply that
     $
     (f(x)-K)(K-x)>(L-\varepsilon)(K-x)^2>0$, $x\in (K-\theta, K)\cup  (K, K+\theta),
     $
  and that conditions \eqref{cond:flhdif}    hold with $a_1=L_0-\varepsilon=a_2$.

 Assume that \eqref{cond:lambda11} holds, set
     \[
\theta_1:=\frac {L_0}{L_0+1}-\alpha-\ell>0
\]
and find $\varepsilon_1\in (0, 1)$ s.t., for each $\varepsilon\in (0, \varepsilon_1)$,
 \[
 \frac {L_0+\varepsilon}{L_0+\varepsilon+1}-\alpha-\ell>\frac{\theta_1}2 \, .
  \]
 Denote, for simplicity of calculations, 
   \[
   \mathcal M:=(1-\alpha-\ell \xi)L_0-\alpha-\ell \xi, \quad \mathcal M(\varepsilon):=(1-\alpha-\ell \xi)(L_0+\varepsilon)-\alpha-\ell \xi,
       \]  
       then, for $\varepsilon\in [0, \varepsilon_1)$, 
       \begin{equation}
  \label{est:mathcalVbel}
  \begin{split}
 \mathcal M(\varepsilon) & >(L_0+\varepsilon)- (\alpha+\ell)(L_0+\varepsilon+1)\\
 &>(L_0+\varepsilon)- \left[ \frac {L_0+\varepsilon}{L_0+\varepsilon+1}-\frac{\theta_1}2  \right](L_0+\varepsilon+1)
 =\frac{\theta_1}2(L_0+\varepsilon+1)>0.
 \end{split}
 \end{equation}
 Acting as in the proof of Theorem \ref{thm:locglgen}, (i),  we obtain
 $
  (G(\alpha+\ell\xi, x)-K)(K-x)>0$, $x\in (K-\theta, K)\cup  (K, K+\theta).$
  Using \eqref{cond:Lipscheps}, for   $ x\in (K-\theta, K)$  we get  $G(\alpha+\ell \xi, x)-K>0$ and 
\[
G(\alpha+\ell \xi, x)-K=(1-\alpha-\ell \xi)(f(x)-K)+(\alpha+\ell \xi)(x-K)\le \mathcal M(\varepsilon) (K-x),
\]
while for   $ x\in (K, K+\theta)$ we obtain $G(\alpha+\ell \xi, x)-K<0$ and 
\[
K-G(\alpha+\ell \xi, x)=(1-\alpha-\ell \xi)(K-f(x))-(\alpha+\ell \xi)(x-K)\le \mathcal M(\varepsilon) (x-K),
\] 
which leads to  $|G(\alpha+\ell \xi, x)-K|<\mathcal M(\varepsilon) |K-x|.$ Now, 
\begin{eqnarray*}
\ln \mathcal M(\varepsilon)=\ln \left[ \mathcal M \times\left(1+\frac{(1-\alpha-\ell \xi)\varepsilon}{ \mathcal M}   \right)\right]=\ln \mathcal M+
\ln  \left(1+\frac{(1-\alpha-\ell \xi)\varepsilon}{ \mathcal M}   \right). 
   \end{eqnarray*}
Choosing $\varepsilon<\frac{\theta_1\lambda_1}4(L_0+1)$, where $\lambda_1$ is from \eqref{cond:lambda11}, applying the inequality $\ln (1+x)<x$, $|x|<1$ and \eqref{est:mathcalVbel},  we arrive at
\[
\frac{(1-\alpha-\ell \xi)\varepsilon}{ \mathcal M}   \le \frac{\varepsilon}{\frac{\theta_1}2(L_0+1) }<\frac{\lambda_1}2<1, \quad \ln  \left(1+\frac{(1-\alpha-\ell \xi)\varepsilon}{ \mathcal M}   \right)<\frac{(1-\alpha-\ell \xi)\varepsilon}{ \mathcal M} <\frac{\lambda_1}2,
\]
and then, using \eqref{cond:lambda11}, we get 
$\displaystyle
\mathbb E \ln \mathcal M(\varepsilon)<\mathbb E\ln \mathcal M+
 \ln  \left(1+\frac{(1-\alpha-\ell \xi)\varepsilon}{ \mathcal M}   \right)\le -\frac{\lambda_1}2.
$

When we do not  assume that $\alpha+\ell< \frac {L_0}{L_0+1}$ and only \eqref{cond:lambda1} holds, we cannot guarantee that $\mathcal M(\varepsilon)>0$. In this case, by \eqref{cond:Lipscheps}, for   $ x\in (K-\theta, K)$  and $G(\alpha+\ell \xi, x)-K>0$, we get
\begin{equation*}
|G(\alpha+\ell \xi, x)-K|=(1-\alpha-\ell \xi)(f(x)-K)+(\alpha+\ell \xi)(x-K)\le (1-\alpha-\ell \xi)(L_0+\varepsilon)(K-x),
\end{equation*}
while, for   $ x\in (K-\theta, K)$  and $G(\alpha+\ell \xi, x)-K<0$, by  \eqref{cond:Lipscheps} we get $f(x)>K$, and then
\[
|G(\alpha+\ell \xi, x)-K|=-(1-\alpha-\ell \xi)(f(x)-K)-(\alpha+\ell \xi)(x-K)\le (\alpha+\ell \xi)(K-x).
\]
Similar estimates  are applied to two other cases: $ x\in (K, K+\theta)$,  $G(\alpha+\ell \xi, x)-K>0$ and $G(\alpha+\ell \xi, x)-K<0$. All the above gives us 
    $ |G(\alpha+\ell \xi, x)-K|<\max\left\{(1-\alpha-\ell \xi)(L_0+\varepsilon),\, (\alpha+\ell \xi)\right\}|K-x|$.
Assuming that $\varepsilon<\frac {\lambda_1}2(L_0)$ and using $\ln(1+x)<x$ for $x \in (0,1)$,  we obtain
\[
\ln \left[(1-\alpha-\ell \xi)(L_0+\varepsilon)\right]=\ln \left[(1-\alpha-\ell \xi) L_0\right]+\ln \left[1+\frac{\varepsilon}{L_0}  \right]<\ln \left[(1-\alpha-\ell \xi)L_0\right]+\lambda_1/2.
\]

Applying the inequality $\max\{a+\epsilon, b\}\le \max\{a, b\}+\epsilon$, $\epsilon>0$, and \eqref{cond:lambda1}, we conclude
\begin{equation*}
\begin{split}
&\mathbb E\ln \max\left\{(1-\alpha-\ell \xi)(L_0+\varepsilon),\, (\alpha+\ell \xi)\right\}\le \mathbb E\max\left\{\ln \left[(1-\alpha-\ell \xi)L_0\right]+\lambda_1/2,\, \ln  (\alpha+\ell \xi)\right\}\\
\le &\mathbb E \max\left\{\ln \left[(1-\alpha-\ell \xi)L_0\right],\, \ln  (\alpha+\ell \xi)\right\} +\lambda_1/2<-\lambda_1/2 .
\end{split}
\end{equation*}
In both cases the rest of the proof  is the same as in Theorem~\ref{thm:locglgen}. 

\subsection{Proof of Theorem \ref{thm:contdiscrdistr}}
Denote by $\psi$ the probability density function (or the probability mass function in a discrete case) of the random variable $\xi$ and let $\mu_2$  be its second moment. Since the distribution is symmetric, we have $\psi(u)=\psi(-u)$, $u\in [-1,1]$, 
 $\mu_2=\int_{-1}^1u^2\psi(u)du$ in the continuous case and $\mu_2=\sum_{i=1}^\infty u_i^2\psi(u_i)$ in the discrete case.
  Choose 
\begin{eqnarray}
\label{def:0ellal}
&&0<\ell_0<\min\left\{\frac 2{\mu_2}, \,\, \frac 1{L_0+1},   \,\, \frac {L_0-1}{(1+\mu_2/2)(L_0+1)} \right\}, \\
\label{def:0ellal1}&&\alpha\in \left(\alpha_0-\frac{\ell_0^2\mu_2}2, \,\, \alpha_0  \right), 
\quad  \ell \in \left(\ell_0, \,\,\min\left\{\alpha, \frac 1{L_0+1}\right\}\right).
\end{eqnarray}
Note  that the second  interval in \eqref{def:0ellal1} is not empty. Indeed, 
\[
\alpha_0-\frac{\ell_0^2\mu_2}2>\ell_0, \quad \mbox{since} \quad \ell_0\left[1+ \frac{\ell_0\mu_2}2 \right]<\ell_0\left[1+ \frac{\mu_2}2 \right]< \frac {L_0-1}{(1+\mu_2/2)(L_0+1)}\left[1+ \frac{\mu_2}2 \right]=\alpha_0,\]
$\alpha>\alpha_0-\frac{\ell_0^2\mu_2}2>\ell_0$ and $\ell_0<\frac 1{L_0+1}$. Also,
$
\alpha+\ell<\alpha_0+\ell<\frac{L_0}{L_0+1},\quad \alpha+\ell>\alpha_0-\frac{\ell_0^2\mu_2}2+\ell_0>\alpha_0,
$
where the second inequality  is true since  $\ell_0<\frac 2{\mu_2}$. 
So we need to prove only the second relation in \eqref{cond:lambda11}.

By Lemma~\ref{lem:mathcalL}~(vi), we have $\mathcal L_0(\alpha_0)
=1$, so $\ln \mathcal L_0(\alpha_0)=0$. For any $\alpha, \ell$ satisfying  \eqref{def:0ellal1}, we get $\alpha+\ell<\frac{L_0}{L_0+1}$, \,
  $\mathcal L_0(\alpha+\ell \xi)
  >0$,\,
$
\mathcal L_0(\alpha+\ell \xi)=\mathcal L_0(\alpha)-\ell (L_0+1)\xi$, \,
 $\mathcal L_0(\alpha_0+\ell \xi)=1-\ell (L_0+1)\xi$, \, $\mathcal L_0^2(\alpha)>\mathcal L_0^2(\alpha_0)=1,
$
\begin{equation}
\begin{split}
\label{est:mathcalL0}
0 & <\mathcal L_0^2(\alpha)-1=[L_0-\alpha(1+L_0)+1][L_0-\alpha(1+L_0)-1]=[1-\alpha](1+L_0)\left[\frac{L_0-1}{(1+L_0)}-\alpha\right](1+L_0)\\&=(1+L_0)^2(1-\alpha)(\alpha_0-\alpha)
<(1+L_0)^2(\alpha_0-\alpha).
 \end{split}
  \end{equation}
Applying the inequality $\ln (1-x)<-x$, $x\in (0, 1)$,  we get, for $u\in [-1, 1]$, 
\begin{equation}
\begin{split}
\label{calc:Lln}
\ln [1-\ell^2 (L_0+1)^2u^2]&<- \ell^2 (L_0+1)^2u^2,\\
 \ln [\mathcal L_0^2(\alpha)-\ell^2 (L_0+1)^2u^2]&<\ln [1+(L_0+1)^2(\alpha_0-\alpha)-\ell^2 (L_0+1)^2u^2]\\&<(L_0+1)^2(\alpha_0-\alpha)-\ell^2 (L_0+1)^2u^2.
\end{split}
\end{equation}
Let now $\xi$ have a continuous distribution, then $\mu_2=\int_{-1}^1u^2\psi(u)du$. Applying \eqref{calc:Lln}, we get 
\begin{equation}
\begin{split}
\label{calc:L0ln}
&\mathbb E \ln [\mathcal L_0(\alpha+\ell\xi)]=\mathbb E \ln [L_0-(\alpha+\ell\xi)(L_0+1)]
=\int_{-1}^1\ln [\mathcal L_0(\alpha)-\ell (L_0+1)u]\psi(u)du\\
&=\int_{0}^1\ln [\mathcal L_0(\alpha)-\ell (L_0+1)u]\psi(u)du-\int_{1}^0\ln [\mathcal L_0(\alpha)+\ell (L_0+1)y]\psi(-y)dy\\
& = \int_{0}^1\ln [\mathcal L_0^2(\alpha)-\ell^2 (L_0+1)^2u^2]\psi(u)du< \int_{0}^1\ln [1+(L_0+1)^2(\alpha_0-\alpha)-\ell^2 (L_0+1)^2u^2]\psi(u)du\\
 &<(L_0+1)^2(\alpha_0-\alpha)\int_{0}^1\psi(u)du-\ell^2_0 (L_0+1)^2\int_{0}^1u^2\psi(u)du=\frac{(L_0+1)^2(\alpha_0-\alpha)}{2}- \frac{\ell^2_0 (L_0+1)^2\mu_2}2\\
 &<\frac{(L_0+1)^2\frac{\ell_0^2\mu_2}2}{2}- \frac{\ell^2_0 (L_0+1)^2\mu_2}2= - \frac{\ell^2_0 (L_0+1)^2\mu_2}4<0,
\end{split}
\end{equation}
which proves the second inequality in \eqref{cond:lambda11}.

Let now $\xi$ be a discrete random variable  with an at most countable  number of states 
$$\{u_1,  -u_1, \dots, u_m, -u_m, \dots\}, \quad u_i\in [-1,1].$$ 
Recall that its probability mass  function $\psi(u)$, $u\in \mathbb R$, is defined as follows: $\psi(u)=0$ when $u\neq u_i$, $\psi(\pm u_i)=P\{\xi=\pm u_i\}$, where $\displaystyle \sum_{m=1}^{\infty} \psi(u_m)= \frac{1}{2}$.
Choose  $\alpha, \ell$ as in \eqref{def:0ellal} and \eqref{def:0ellal1} and 
denote $\displaystyle H:=\max_{u\in [-1, 1]}\left|\ln [\mathcal L_0(\alpha)-\ell (L_0+1)u]\right|$.
Since $\displaystyle \sum_{i=1}^\infty \psi(u_i)$ is convergent, we can find $N_1\in\mathbb N$ s.t. 
$
\displaystyle \sum_{i=N_1+1}^\infty \psi(u_i)<\ell_0^2(L_0+1)^2\mu_2/(16H), \quad \sum_{i=1}^{N_1} u_i^2\psi(u_i)>7\mu_2/16.
$

The series $\displaystyle \sum_{i=1}^\infty\ln [\mathcal L_0(\alpha)\pm \ell (L_0+1)u_i]\psi (u_i)$ is absolutely 
convergent, so we can estimate
$$
\sum_{i=N_1+1}^\infty\ln [\mathcal L_0(\alpha) \pm \ell (L_0+1)u_i]\psi (u_i)
\le H\sum_{i=N_1+1}^\infty\psi (u_i)<\ell_0^2(L_0+1)^2\mu_2/16.$$
Further,
\begin{eqnarray*}
\mathbb E \ln [\mathcal L_0(\alpha+\ell\xi)]&=&\sum_{i=1}^{N_1}\ln [\mathcal L_0(\alpha)-\ell (L_0+1)u_i]\psi (u_i)+\sum_{i=1}^{N_1}\ln [\mathcal L_0(\alpha)+\ell (L_0+1)u_i]\psi (u_i)\\
&&+\sum_{i=N_1+1}^\infty\ln [\mathcal L_0(\alpha)-\ell (L_0+1)u_i]\psi (u_i)+\sum_{N_1+1}^\infty\ln [\mathcal L_0(\alpha)+\ell (L_0+1)u_i]\psi (u_i)\\
& < & \sum_{i=1}^{N_1}\ln [\mathcal L_0^2(\alpha)-\ell ^2(L_0+1)^2u_i]^2\psi (u_i)+\ell_0^2(L_0+1)^2\mu_2/8.
\end{eqnarray*}
Now, applying \eqref{def:0ellal}, \eqref{def:0ellal1}, \eqref{est:mathcalL0} and acting as in \eqref{calc:L0ln}, we arrive at 
\begin{equation*}
\begin{split}
\mathbb E \ln [\mathcal L_0(\alpha+\ell\xi)] 
& <\sum_{i=1}^{N_1}\ln [1+(L_0+1)^2(\alpha_0-\alpha)-\ell^2 (L_0+1)^2 u_i^2]\psi (u_i)+
\ell_0^2(L_0+1)^2\mu_2/8\\
 &<(L_0+1)^2(\alpha_0-\alpha)\sum_{i=1}^{N_1}\psi (u_i) -\ell^2_0 (L_0+1)^2\sum_{i=1}^{N_1} u_i^2\psi (u_i)+\ell_0^2(L_0+1)^2\mu_2/8\\
 &
 \le \frac{(L_0+1)^2}{2}\left[\ell_0^2\mu_2/2-\ell^2_0 7\mu_2/8 +\ell^2_0 \mu_2/4 \right]= -\ell_0^2(L_0+1)^2\mu_2/8 <0, 
   \end{split}
  \end{equation*}
which is the second inequality in \eqref{cond:lambda11}. The reference to Theorem~\ref{thm:Singer} concludes the proof.

\subsection{Proof of Proposition~\ref{prop:Globdet}}
Let and $\mathcal L^{\pm}$ be defined by \eqref{def:mathcalL_+12} and  $\Psi (\cdot, \cdot)$ 
by \eqref{def:Psi}.
Since $\underline \beta_0=\Psi (L^-, L^+)$, see  \eqref{def:beta0}, by Lemma~\ref{lem:mathcalL}~(vi), we have $\mathcal L^{-}(\underline \beta_0)\mathcal L^{+}(\underline \beta_0)=1$,  and also $\mathcal L^{\pm}(\alpha)<\mathcal L^{\pm}(\underline \beta_0)$.  Define 
\begin{equation*}
\begin{split}
\phi(x)=\mathcal L^{-}(\underline \beta_0)(K-x)+K, \, x_m \le x\le K, \,\, \phi(x)=-\mathcal L^{+}(\underline \beta_0)(x-K)+K, \,   K\le x\le \mathcal L^{-}(\underline \beta_0)(K-x_m)+K,
\end{split}
\end{equation*}
which is decreasing, $\phi(x)>K$, $ x_m \le  x < K$, $\phi(x)<K$, $ K < x \le \mathcal L^{-}(\underline \beta_0)(K-x_m)+K$
 and $\phi(\phi(x))=x$.
 Since $\phi(\mathcal L^{-}(\underline \beta_0)(K-x_m)+K)=x_m>0$, the function $\phi$ is also positive.  
By \eqref{cond:Lglob2} and since  $\mathcal L^{\pm}(\alpha)<\mathcal L^{\pm}(\underline \beta_0)$, we get, for $\alpha>\underline \beta_0$,  
$
G(\alpha, x)\le \mathcal L^{-}(\alpha) (K-x)+K<\phi(x)$, $ x\in (0, K)$, $G(\alpha, x)\ge- \mathcal L^{+}(\alpha) (x-K)+K>\phi(x)$, $x\in (K, f_m),
$
and  the result follows from Lemma~\ref{lem:Cull2}.


\subsection{Proof of Proposition~\ref{prop:envG}}
To construct an envelope $\phi$ for $G(\alpha, x)$, $\alpha > \alpha_0$, we estimate 
\[
G(\alpha_0, x)-G(\alpha_0, a_i)\le L_i^{-}(1-\alpha_0)(a_i-x)-\alpha_0(a_i-x)=\mathcal L^{-}_i(\alpha_0)(a_i-x),\quad a_{i+1}\le x\le a_i, 
\] 
\[
G(\alpha_0, b_i(\alpha_0))-G(\alpha_0, x)\le L_i^{+}(1-\alpha)(x-b_i(\alpha_0))-\alpha(x-b_i(\alpha_0))=\mathcal L^{+}_i(\alpha_0)(x-b_i(\alpha_0)),\quad b_{i}(\alpha_0)\le x. 
\]
Set  $\alpha_i :=\Psi(L_i^-, L_i^+)$, then \eqref{cond:alpha0gl} implies that  $\alpha_i\le \alpha_0$.
Since $\mathcal L^{\pm}_i(\alpha)\ge  \mathcal L^{\pm }_i(\alpha_0)$ when $\alpha \le \alpha_0$,  and $ \mathcal L^{-}_i(\alpha_i) \mathcal L^{+}_i(\alpha_i)=1$, see Lemma~\ref{lem:mathcalL}~(vi), 
we  get that $ \mathcal L^{-}_i(\alpha_0) \mathcal L^{+}_i(\alpha_0)\le 1$. Setting
$
C_i^-:=\mathcal L^{-}_i(\alpha_0)$, $ C_i^+:=\left[\mathcal L^{-}_i(\alpha_0)\right]^{-1}, 
$
we get $C_i^-C_i^+=1$, $C_i^+\ge \mathcal L^{+}_i(\alpha_0)$ and 
\[
G(\alpha_0, x)-G(\alpha_0, a_i)\le C_i^- (a_i-x),\,\, a_{i+1}\le x\le a_i, \quad G(\alpha_0, b_i(\alpha_0))-G(\alpha_0, x)\le C_i^+(x-b_i(\alpha_0)), \,  b_{i}(\alpha_0)\le x.
\]
Define $\phi$, as in \eqref{def:phi}, and by straightforward calculations, show that 
$
G(\alpha_0, x)\le \phi(x)$, $x\in (0, K)$, $G(\alpha_0, x)\ge \phi(x)$, $x>K. 
$
Since for   $\alpha > \alpha_0$,
$
G(\alpha, x)< G(\alpha_0, x), \, x\in (0, K), \quad G(\alpha_0, x) <G(\alpha, x), \,  x>K, 
$
we conclude that  $\phi$  is an envelope for $G(\alpha, x)$,  therefore $\lim\limits_{n\to \infty}x_n=K$, 
by Lemma~\ref{lem:Cull2}.


\subsection{Remark to Proposition \ref{prop:envG}}
\begin{remark}
\label{rem:glob}
If there is $i_0<m$ s.t.  $\alpha_0<\Psi \left(L_{i_0}^-, L_{i_0}^+ \right)$,  we can find a bigger parameter $\bar \alpha$ which guarantees global stability of the solution to \eqref{eq:PBCintrdet}.  To show that we denote  
\begin{equation*}
\begin{split}
&f(x)-f(y)\le L_i^-(y-x), \,\, a_{i+1}\le x<y\le a_i, \quad i=0, 1, \dots, m-1,\\
&L^+(z), \, z>K, \, \mbox{be a Lipschitz constant for }  f  \mbox{ s.t }\,\, f(y)-f(x)<L^+(z)(x-y), \,\, \forall x>y\ge z,
\end{split}
\end{equation*}
$b_0=K$, $b_1=b_1(\alpha_0):=\max_{x\in [a_1, K]}\{G(\alpha_0, x)\}$, 
$b_i=b_i(\alpha_{i-1})=\max_{x\in [a_i, K]}\{G(\alpha_i, x ) \}$, where $\alpha_i$ are defined inductively: 
\[
\alpha_1:=\max\{\alpha_0, \Psi(L_1^-, L^+(b_1))\}, \dots, \alpha_k:=\max\{\alpha_0, \Psi(L_i^-, L^+(b_i)), i=1, \dots, k\}.
\]
Set $\displaystyle \bar \alpha:= \max_{ i=0, 1, \dots, m-1}\{\alpha_i\}$, $ \tilde{\mathcal L}^{+}_i(\alpha):=(1-\alpha)L^+(b_{i})-\alpha$ and note that $\lim\limits_{n\to \infty} x_n=K$, for any $\alpha>\bar \alpha$ and $x_0\in (a_1, b_1)$. 
We want to get the same for each $x_0>0$.  Fix some $\alpha>\bar \alpha$ and assume the contrary: for some $k<m$ we have stability on $(a_k, b_k)$ but $\bar x:=\inf\left\{x\in (a_{k+1}, a_k): G^2(\alpha, x)>x  \right\}\ge a_{k+1}.$ This implies that 
$G^2(\alpha, \bar x)=\bar x$.  By the inductive assumption for $(a_k, b_k)$ we get  $G(\alpha, \bar x) >b_k$ .  Also, $G(\alpha, \bar  x) <b_{k+1}$ since  $\alpha>\bar \alpha> \alpha_{k+1}$ and therefore $G(\alpha, \bar x)<G(\alpha_{k+1}, \bar  x)\le b_{k+1}$, where the last inequality holds by the definition of $b_{k+1}$.
Choose $\hat  x\in (\bar  x, a_k)$ s.t. $G(\alpha, \hat  x)\in (b_{k}, b_{k+1})$. Then $\hat x>\bar x$, $G^2(\alpha, \hat x)>\hat x$,  and 
$
G(\alpha, \bar x)-G(\alpha, \hat x)\le \mathcal L^{-}_k(\alpha)(\hat x-\bar x).
$
Assuming $G(\alpha, \hat x)\le G(\alpha, \bar x)$ we get 
\begin{eqnarray*}
&\hat x- G^2(\alpha, \bar x)  & <
G^2(\alpha, \hat x)-G^2(\alpha, \bar x)=G\bigl(\alpha, G(\alpha, \hat x)\bigr)-G\bigl(\alpha, G(\alpha, \bar x)\bigr)\\&&\le \tilde{\mathcal L}_k^{+}(\alpha)[G(\alpha, \bar  x)-G(\alpha, \hat x)]\le \tilde{\mathcal L}_k^{+}(\alpha)\mathcal L^{-}_k(\alpha)(\hat x-\bar x)<\hat x-\bar x \implies G^2(\alpha, \bar x)> \bar x,
\end{eqnarray*}
contradicting to  the definition of  $\bar x$.
If however, $G(\alpha, \hat x)> G(\alpha, \bar x)$, we can find $\tilde x\in (\hat x, K)\in  (\bar x, K)$ s.t. $G(\alpha, \bar  x)=G(\alpha, \tilde  x)$. But then $\bar  x=G^2(\alpha, \bar  x)=G^2(\alpha, \tilde  x)>\tilde  x$,
which contradicts to the choice of  $\tilde x$.
\end{remark}

\subsection{Proof of Proposition~\ref{prop:3}}
Define 
\[
\mathcal U(x, \alpha_0):=G^2(\alpha_0, x)-x,  \quad x\in [x_m, K].
\]
Note that $\mathcal U(K, \alpha_0)=0$ and $\mathcal U_x'(x, \alpha_0)=G'\bigl(\alpha_0, G(\alpha_0, x)\bigr)G'(\alpha_0, x)-1=[(1-\alpha_0)f'( G(\alpha_0, x)) +\alpha_0][(1-\alpha_0)f'(x) +\alpha_0]-1.$
Fix some $x\in (x_{\max}, K)$ and note that the  equation $V(u)=0$, with 
\[
V(u)=V(x, u)  :=  [f'( G(\alpha_0, x))+u(1-f'( G(\alpha_0, x)))][f'(x) +u(1-f'(x))]-1
\]
has two real roots, $1$ and $\displaystyle{\Psi \left(f'(G(\alpha_0, x)), f'(x) \right)=\frac{f'(G(\alpha_0, x)) f'(x)-1 }{(1-f'(G(\alpha_0, x)), (1-f'(x)) }} \leq 1$, 
by Lemma~\ref{lem:mathcalL}~(v).
Also, $V(u)<0$ when $u\in \bigl(\Psi \left(f'(G(\alpha_0, x), f'(x) \right), 1\bigr)$. By \eqref{cond:alpha0}, we have 
$\alpha_0\in \left(\Psi \left(f'(G(\alpha_0, x), f'(x)) \right), 1\right)$, which  implies that 
$
\mathcal U'_x(x, \alpha_0)=V(\alpha_0)< 0,
$
so $\mathcal U(x, \alpha_0)$ decreases in $x$. Therefore, for each $x\in (x_{\max}, K)$,  
$
G^2(\alpha_0, x)-x=\mathcal U(x, \alpha_0)> \mathcal U(K, \alpha_0)=0.
$
Set  
\[
G_{\max}(\alpha_0)=\max_{x\in [0, K]}G(\alpha_0,x),  \,  x_{Gmax} \,\, \mbox{is  the largest point of maximum of $G(\alpha_0,\cdot)$ on (0, K)}.
\]

For each  $y\in (K, G_{\max}(\alpha_0))$ there is $x\in (0, K)$ s.t. $y=G(\alpha_0, x)$. Due to continuity  we can choose $x\in (x_{Gmax}, K)\subseteq (x_{\max}, K)$.  Since $G^2(\alpha_0, x)-x=\mathcal U(x, \alpha_0)>0$, we conclude that  $G^2(\alpha_0, x)>x$ and therefore $x<G(\alpha_0, y)$.  If $G^2(\alpha_0, y)>y$,  there is a point $\hat x\in (G(\alpha_0, y), K)$ s.t. $y=G(\alpha_0, x)=G(\alpha_0, \hat x)$, so  $x<G^2(\alpha_0, x)=G^2(\alpha_0, \hat x)=G(\alpha_0, y)<\hat x$, or $G^2(\alpha_0, \hat x)<\hat x$, which is a contradiction to the case proved above.

When $x\in (0, x_{\max})$ there exists $\bar x\in (x_{\max}, K)$ s.t. $G(\alpha_0, x)=G(\alpha_0, \bar x)$ and we are in the first case.  The case $x>G_{\max}(\alpha_0)$ is treated as in Lemma~\ref{lem:aux21}.

Application of Lemma~\ref{lem:Gull} proves that any  control  $\alpha>\alpha_0$ guarantees global stability.


\begin{thebibliography}{99}



\bibitem{BR0}
G. Berkolaiko and A. Rodkina, 
Asymptotic behavior of solutions to linear discrete stochastic equation, in 
Proceedings of {\it the International Conference ``2004-Dynamical Systems and Applications"}, 
Antalya, Turkey, 5-10 July 2004, 614--623.

\bibitem{BR1}
G. Berkolaiko and A. Rodkina,
Almost sure convergence of solutions to non-homogeneous stochastic difference equation, 
{\em J. Difference Equ. Appl.} {\bf 12} (2006),
535--553. 


\bibitem{BKR2016}
E. Braverman, C. Kelly and  A. Rodkina,
Stabilization of difference equations with noisy prediction-based control,
{\it Physica D},  {\bf 326} (2016), 21--31.

\bibitem{BKR2020}
E. Braverman, C. Kelly and A. Rodkina,
Stabilization of cycles with stochastic
prediction-based and target-oriented  control, 
{\em Chaos} {\bf 30} (2020), 15pp.

\bibitem{BL2012}
E.~Braverman and E.~Liz,
\newblock On stabilization of equilibria using predictive control with and
without pulses,
\newblock {\em Comput. Math. Appl.} {\bf 64} (2012), 
2192--2201.

\bibitem{BRAllee}
E. Braverman and A. Rodkina,
Stochastic difference equations with the Allee effect,   
{\em Discrete Contin. Dyn. Syst. Ser. A} {\bf 36}(11) (2016),  5929--5949.

\bibitem{BR2019}
\newblock E. Braverman and A. Rodkina,
\newblock {Stochastic control stabilizing unstable or chaotic maps},
 {\em J. Difference Equ. Appl.} { \bf 25} (2019), 151--178.
 
 
%

\bibitem{BR2022}
E. Braverman and A. Rodkina,
Stabilizing multiple equilibria and cycles with noisy prediction-based control, 
 {\em Discrete Contin. Dyn. Syst. Ser. B} {\bf 27} 
(2022),  5419--5446.


\bibitem{Chagas}
T. P. Chagas, P.-A. Bliman, and K. H. Kienitz,  
Stabilization of periodic orbits of discrete-time dynamical systems using the prediction-based control: new control law and practical aspects, {\em J. Franklin Inst.} {\bf 355} (2018), 
4771--4793. 
 


\bibitem {Elaydi} 
S. Elaydi, An Introduction to Difference Equations, 
{\em Undergraduate Texts in Mathematics}, Springer, New York, 2005. 

 
\bibitem {ElaydiSacker} 
	S. Elaydi and R. Sacker, 
	Basin of attraction of periodic orbits of maps on real line,
   {\em J. Difference Equ. Appl.} {\bf 10} 
	(2004), 881--888.
  
\bibitem{Coppel} 
W.\,A. Coppel,  The solution of equations by
iteration, {\em Proc. Cambridge Philos. Soc.} {\bf 51} (1955), 41--43.	
 
  \bibitem{Gull} 
	P. Cull,
	Global stability of population models, 
{\em	Bull. Math. Biol.} {\bf 43} (1981), 47--58. 
  
	
\bibitem{Cull2} 
P. Cull, 
Population models: stability in one dimension,
{\em Bull. Math. Biol.} {\bf 69} (2007), 989--1017.
  
%

\bibitem{Kh}
R. Has’minski, Stability of Systems of Differential Equations
under Random Perturbations of Their Parameters, Nauka,
Moscow, 1969, 367 pp. 

 \bibitem{Medv}
P. Hitczenko and G. Medvedev, 
Stability of equilibria of randomly perturbed maps,  
{\em Discrete Contin. Dyn. Syst. Ser. B} {\bf 22}(2) (2017), 269--281.

\bibitem{K} 
H. Kesten, Random difference equations and renewal theory for the
product of random matrices, {\em Acta Math.} {\bf 131} (1973), 207--248.


 
 \bibitem{FL2010}
E.~Liz and D.~Franco, 
\newblock Global stabilization of fixed points using predictive control,
\newblock {\em Chaos} {\bf 20} (2010), 023124, 9 pages.

\bibitem{mb} 
J.~G. Milton, J. B\'elair, 
Chaos, noise, and extinction in models of population growth,
{\em Theor. Popul. Biol. } {\bf  37} (1990) 273--290.

\bibitem{polyak}   
B.\,T. Polyak,  Chaos stabilization by predictive control,
{\em Autom. Remote Control}  {\bf 66} (2005), 1791-1804.

\bibitem{Rich}
S. Rich, A. Hutt, F.\,K. Skinner, T.\,A. Valiante, and J. Lefebvre,
Neurostimulation stabilizes spiking neural networks by disrupting seizure-like oscillatory transitions, {\em Sci Rep} {\bf 10} (2020), 15408. 
https://doi.org/10.1038/s41598-020-72335-6

\bibitem{Shark}  
A.\,N. Sharkovsky, Yu.\,L. Maistrenko and E.\,Yu. Romanenko, 
Difference Equations and Their Applications, Kluwer Academic, Dordrecht, 1993.


\bibitem{Shiryaev96} 
A.~N.~Shiryaev.
\newblock \emph{Probability}. 2nd edition,
\newblock Springer, Berlin, 1996.

\bibitem{Singer}
D.~Singer, 
\newblock Stable orbits and bifurcation of maps of the interval,
\newblock {\em SIAM J. Appl. Math.} {\bf 35}(2)  (1978), 260--267.




\bibitem{sousa} 
M. de Sousa Vieira and A.\,J. Lichtenberg,   Controlling chaos using nonlinear feedback with delay,
{\em Phys. Rev. E} {\bf 54} (1996), 1200--1207.

\bibitem{uy99} 
T.~Ushio and S.~Yamamoto,   
\newblock Prediction-based control of chaos, 
\newblock {\em Phys. Lett. A} {\bf  264} (1999), 30--35.
\end{thebibliography}
\end{document}